\documentclass{amsart}
\usepackage{a4wide}
\usepackage[latin1]{inputenc}
\usepackage{graphicx,epsfig,subfig}

\newtheorem{theorem}{Theorem}

\newtheorem{ncor}{Corollary}
\newtheorem{lemma}{Lemma}
\newtheorem{prop}{Proposition}

\theoremstyle{remark}

\newtheorem{rem}{Remark\!\!}

\newtheorem{ack}{Acknowledgements\!\!}

\renewcommand{\(}{\left(}
\renewcommand{\)}{\right)}

\newcommand{\C}{$C[0,\infty)$} 
\newcommand{\p}[1]{{\bf P}\left\{#1\right\}} 
\newcommand{\Ord}[1]{{O}\left(#1\right)}  
  
\newcommand{\rcp}[1]{\frac{1}{#1}} 
\newcommand{\E}{{\bf E}} 

\newcommand{\bth}{\begin{theorem}} 
\newcommand{\eth}{\end{theorem}} 
\newcommand{\bl}{\begin{lemma}} 
\newcommand{\el}{\end{lemma}}
\newcommand{\bpr}{\begin{prop}}
\newcommand{\epr}{\end{prop}}
\newcommand{\bcor}{\begin{ncor}}
\newcommand{\ecor}{\end{ncor}}
\newcommand{\brem}{\begin{rem}} 
\newcommand{\erem}{\end{rem}} 
\newcommand{\bpf}{\medskip\noindent{\it Proof:}\ \ } 
\newcommand{\epf}{\hfill{ }$\Box$\vspace*{1.5ex}} 
\newcommand{\eps}{\varepsilon}

\newcommand{\cw}{\stackrel{w}{\longrightarrow}} 
\newcommand{\pdiff}[3]{\frac{\partial^{#3}#1}{\partial #2^{#3}}} 
\newcommand{\Tn}{${\mathcal Y}_n$}
 
\newcommand{\nti}{{n\to\infty}} 
\newcommand{\U}[1]{\left[ #1 \right]_{u=1}}
\newcommand{\btr}[1]{\left| #1 \right|}

\author{Michael Drmota and Bernhard Gittenberger}
\title{The Shape of Unlabeled Rooted Random Trees}
\date{\today}
\thanks{Address: Institut für Diskrete Mathematik und Geometrie,
Technische Universität Wien, Wiedner Hauptstr. 8-10/104, A-1040 Wien,
Austria. \\
Email: {\tt michael.drmota@tuwien.ac.at, gittenberger@dmg.tuwien.ac.at}}
\thanks{This work has been supported by the Austrian Science Foundation FWF, NFN-Project S9604, 
as well as by the EU FP6-NEST-Adventure Programme, contract no. 028875 (NEST)}


\begin{document}

\begin{abstract}
We consider the number of nodes in the levels of unlabelled rooted random trees and show that the
stochastic process given by the properly scaled level sizes weakly converges to the local time of
a standard Brownian excursion. Furthermore we compute the average and the distribution of the 
height of such trees. These results extend existing results for conditioned Galton-Watson trees
and forests to the case of unlabelled rooted trees and show that they behave in this respect 
essentially like a conditioned Galton-Watson process. 
\end{abstract}

\maketitle

\section{Introduction}

We consider the profile and height of unlabelled rooted random trees. This kind of trees is 
also called P\'olya trees, because the enumeration theory developed by P\'olya allows an
analytical treatment of this class of trees by means of generating functions (see \cite{Pol37}).
The profile of a rooted tree $T$ is defined as follows. First we define the $k$-th level of $T$ to
be the set of all nodes having distance $k$ from the root (where we use the usual shortest path
graph metric). Let $L_k(T)$ denote the number of nodes of the $k$-th level. The \emph{profile} of
$T$ is the sequence $(L_k(T))_{k\ge 0}$. For a random tree this sequence becomes a stochastic
process.

The first investigations of the profile of random trees seem to go back to
Stepanov \cite{Stepanov69} who derived explicit formulas for the distribution of the
size of one level. Further papers deal mainly with 
simply generated trees as defined by Meir and Moon \cite{MM78}. Note that simply generated trees
are defined by the functional equation 
\begin{equation} \label{simplygen}
y(x)=x\phi(y(x))
\end{equation}
for their generating function but can also be viewed as family
trees of a Galton-Watson process conditioned on the total progeny. Kolchin (see
\cite{Kol77,Kol86}) related the level size distributions to distributions
occurring in particle allocation schemes. Later Tak\'acs \cite{Ta91} derived another
expression for the level sizes by means of generating functions. 
Aldous \cite{Al91} conjectured two functional limit theorems for the profile in two
different ranges which were proved in \cite{DG97,G98}. The first author \cite{Dr96}
studied restrictions of the profile to nodes of fixed degree. An extension to
random forests of simply generated trees is given by the second author \cite{G02}. 

Later other tree classes have been considered as well. The profile of random
binary search trees has been first studied by Chauvin et al. \cite{CDJ01} and later by Drmota and
Hwang (see \cite{Dr04} and \cite{DH05a}). Random recursive trees have been investigated
recently by Drmota and Hwang \cite{DH05b} and van der Hofstad et al. \cite{HHM02}. 
Related research was done by Chauvin et al. \cite{CDJ01}, Fuchs et al. \cite{FHN06}, Hwang
\cite{Hw05, Hw07}, Louchard et al. \cite{LST}, and Nicod\`eme \cite{Ni05}. 
Extremal studies of the profile (called the width of trees)
of simply generated trees have been started by Odlyzko and Wilf \cite{OW87}. 
The distribution including moment
convergence has been presented independently in Chassaing et al. \cite{CMY} and the authors of
this paper \cite{DG04}. For other tree classes we refer to the work of Devroye and Hwang
\cite{DeHw06} and Drmota and Hwang \cite{DH05b}. A general overview on random trees which also
strongly highlights the profile of trees can be found in the first author's book on random trees
\cite{Drm09}.  

Whereas simply generated trees have an average height of order $\sqrt n$, the other tree classes
mentioned above have height of order $\log n$. P\'olya trees do 
not belong to the class of simply generated trees. Since we are not aware of any rigorous proof of
this assertion in the literature, we will present a (rather simple) 
proof of this in the next section. 
To our knowledge, so far the fact was only underpinned by the following 
argument concerning the generating functions of both tree classes. 
The argument works for many tree classes considered in the literature, e.g. Cayley trees, 
plane trees, Motzkin trees, binary trees, and many more. In all these cases the 
function $\phi(y)$ in \eqref{simplygen} is entire function or meromorphic. But it is not at all
clear that this argument is still true if $\phi(y)$ has a more complicated singularity structure.
In these cases (i.e., entire or meromorphic $\phi(y)$), 
the generating functions enumerating the number of simply generated trees and 
P\'olya trees, respectively, have a fundamentally different singularity
structure. Whereas the first one has one or a finite number of singularities (the latter occurs in
the periodic case) on the circle of convergence and allows analytic continuation to a slit plane
(with the possible exception of finitely many isolated singularities which are of algebraic type
even if the function itself is not algebraic), 
the generating function associated to P\'olya trees is much more
complicated. In fact, for the latter function the unit circle is a natural boundary (i.e., no
analytic continuation beyond it is possible). There is exactly one singularity on the circle of
convergence of the power series expansion at 0, but the analytic continuation has an infinite
number of singularities inside the unit circle. Each point on the unit circle is an accumulation
point of the set of singularities. These facts follow from the functional equation defining this
generating function and the fact that the power series expansion around zero has radius of
convergence strictly smaller than one (see next section). It also involves an analytically
complicated structure like the cycle index of the symmetric group. Due to this difference with
respect to the analytic behaviour of the generating function P\'olya trees are not
simply generated and therefore they cannot be represented as branching processes. 

Note that the rather complicated singularity structure does not affect the asymptotics of
statistical parameters like the number of trees, the profile or the height. The behaviour of these
parameters is determined by the dominant singularity which is, as all other singularities as well,
an isolated singularity, exactly as in the case of simply generated trees. 
Indeed, P\'olya trees behave in many respects similar to simply generated trees (compare with
\cite{RS75, HRS75, MM91,MM92, DG3, Gi06}) or the recent work of Marckert and Miermont
\cite{MaMi09} who showed that binary unlabelled trees converge in some sense to the continuum 
random tree, i.e., the same limit as that of simply generated trees. Moreover, Broutin and Flajolet
\cite{broutinflajolet} showed $\sqrt n$-behaviour for the height of binary unlabelled trees. 
Hence it is expected that the order of the height is
$\sqrt n$ as well. In this paper we will give an affirmative answer to this question. This
justifies the choice of $\sqrt n$ for the scaling of the level sizes in the subsequent theorems. 

The plan of the paper is as follows. In the next section we present our main results. Then we will
set up the generating functions for our counting problem of trees with nodes in certain levels
marked. This function is given as solution of a recurrence relation which has to be analyzed in
detail. Knowing the singular behaviour of the considered generating functions allows us to show
that the finite dimensional distributions (fdd's) of the profile, i.e., the distributions of the
sizes of several levels considered simultaneously, converge to the fdd's of Brownian excursion
local time. The singularity analysis is carried out in Section~\ref{mainrec} and the computation
of the fdd's in Section~\ref{fdds}. 
In order to complete the functional limit theorem we need to prove tightness. This
means, roughly speaking, that the sample paths of the process do not have too strong fluctuations
(see \cite{Bi68} for the general theory). 

In the final section we turn to the height. The pioneering work on this topic was done by 
Flajolet and Odlyzko \cite{FO82} and Flajolet et al. \cite{FGOR93} in their studies of 
simply generated trees
where they completed the program started in \cite{BKR}. What we have to do is to show
that the generating function appearing in the analysis of the height has a local structure which
is amenable to the steps carried out in \cite{FO82} and \cite{FGOR93}. 
This is done in the last section and 
leads to average and distribution of the height.

\section{Preliminaries and Results}

First we collect some results for unlabelled unrooted trees. 
Let \Tn\  denote the set of unlabelled rooted trees consisting of $n$ vertices and $y_n$ be the
cardinality of this set. P\'olya \cite{Pol37} already discussed the generating
function  
\[
y(x)=\sum_{n\ge1}y_n x^n
\]
and showed that the radius of convergence $\rho$ satisfies $0<\rho<1$ and
that $x=\rho$ is the only singularity on the circle of convergence
$|z|=\rho$. He also showed that $y(x)$ satisfies the functional equation 
\begin{equation} \label{funeq}
y(x)=x\exp\(\sum_{i\ge 1}\frac{y(x^i)}{i}\).
\end{equation} 
Nowadays, this functional equation is easily derived by using the theory of combinatorial
constructions which is presented in the comprehensive book of Flajolet and Sedgewick
\cite{FlSe10}. Indeed, a P\'olya tree can be viewed as a root with a multiset of P\'olya trees
attached to it. Then the functional equation~\eqref{funeq} pops out immediately from the multiset
construction and its generating function. This functional equation can be used to compute the
coefficients: 
\begin{equation} \label{Reihe}
y(x)=x+x^2+2x^3+4x^4+9x^5+20x^6+48x^7+115x^8+286x^9+719x^{10}+\dots
\end{equation}  

Later Otter \cite{Otter} showed that $y(\rho)=1$ as well as the
asymptotic expansion 
\begin{equation} \label{y-asym}
y(x)=1-b(\rho-x)^{1/2}+c(\rho-x)+d(\rho-x)^{3/2}+\cdots
\end{equation} 
which he used to deduce that
\begin{equation} \label{yn}
y_n\sim\frac{b\sqrt\rho}{2\sqrt\pi}n^{-3/2}\rho^{-n}.
\end{equation} 
Furthermore he calculated the first constants appearing in this expansion:
$\rho\approx0.3383219$, $b\approx2.6811266$, and $c=b^2/3\approx2.3961466$.

We will return to the function $y(x)$ in Section~\ref{mainrec} and list a couple of useful
properties in Lemma~\ref{lem_2} after introducing some notations.

\bth
P\'olya trees are not simply generated. 
\eth

\bpf
Let us assume that P\'olya tree are simply generated. Then the generating function $y(x)$ given by
\eqref{funeq} must have a representation in the form \eqref{simplygen} where $\phi(y)$ is a power
series with non-negative coefficients. 
By \eqref{Reihe} the functional inverse $y^{-1}(x)$ exists and 
we have $y(x)\sim x$ and $y^{-1}(x)\sim x$, as $x\to 0$. This implies $\phi(0)=1$. Plugging 
$y^{-1}(x)$ into \eqref{simplygen} we obtain 
$$
x=y(y^{-1}(x))=y^{-1}(x) \phi(y(y^{-1}(x)))=y^{-1}(x) \phi(x)
$$
and consequently
\begin{align*} 
\phi(x)&=\frac{x}{y^{-1}(x)}=\frac{x}{x-x^2+x^4-x^5+x^6-4x^7+11x^8-18x^9+18x^{10}+\dots} \\
&=1+x+x^2-x^5+3x^6-5x^7+7x^8-8x^9+x^{10}+\dots
\end{align*} 
which violates the requirement of non-negative coefficients for $\phi(x)$.
\epf

\brem
The sequence of the coefficients of $y(x)$ is A000081, that of $y^{-1}(x)$ is A050395 in Sloane's
On-line Encyclopedia of Integer Sequences \cite{Sloane00}
\erem

The height of a tree is the maximal number of edges on a path from the root to another vertex of
the tree. It turns out that the average height is of order $\sqrt n$. 

\bth\label{height}
Let $H_n$ denote the height of an unlabelled rooted random tree with $n$ vertices. Then we
have
\begin{equation} \label{av_height}
\E H_n\sim \frac{2\sqrt\pi}{b\sqrt\rho} \sqrt n. 
\end{equation} 
and 
\begin{equation} \label{dist_height}
\E H_n^r\sim \left( \frac2{b\sqrt\rho}\right)^r
 r(r-1) \Gamma(r/2) \zeta(r)\,  n^{r/2}
\end{equation} 
for every integer $r\ge 2$.
\eth

The proof of this theorem is deferred to the last section, since the proofs of the 
auxiliary lemmas which will eventually establish the assertion will utilize similar techniques as
needed to prove the next three theorems. 

\brem
Note that more information on the limiting distribution is available. Indeed, a local limit
theorem holds as well. Let $y_n^{(h)}$ denote the number of
unlabelled rooted trees with $n$ vertices and height equal to $h$ and let $\delta>0$ 
arbitrary but fixed. If we set $\beta=2\sqrt n/hb\sqrt\rho$, then,
as $\nti$, we have 
$$
\p{H_n=h}=\frac{y_n^{(h)}}{y_n}
\sim 4b\sqrt\frac{\rho\pi^5}{n} \beta^4\sum_{m\ge 1} m^2(2m^2\pi^2\beta^2-3)
e^{-m^2\pi^2\beta^2}
$$
uniformly for $\rcp{\delta\sqrt{\log n}}\le \frac{h}{\sqrt n} \le \delta\sqrt{\log n}$.
A rigorous proof of this theorem was given by Broutin and Flajolet \cite{broutinflajolet}
for binary unlabelled trees. They also provide a moment convergence theorem, a weak limit theorem
as well as large deviation results. 
\erem

Let $L_n(t)$ denote the number of nodes at distance $t$ from the root
of a randomly chosen unlabelled rooted tree of size $n$. If $t$ is not an integer, then
define $L_n(t)$ by linear interpolation: 
\begin{equation} \label{L_n}
L_n(t)=(\lfloor t\rfloor+1-t)L_n(\lfloor t\rfloor)+(t-\lfloor t\rfloor) 
L_n(\lfloor t\rfloor+1), \quad t\ge0.
\end{equation} 

We will show the following theorem.

\bth\label{theo_1}
Let 
$$
l_n(t)=\rcp{\sqrt n}L_n\(t\sqrt n\).
$$ 
Then $l_n(t)$ satisfies the following functional limit theorem:  
\[
(l_n(t))_{t\ge 0} \cw \(\frac {b\sqrt \rho}{2\sqrt2} \cdot l\left(\frac {b\sqrt\rho}{2\sqrt2}\cdot
t\right) \)_{t\ge 0}
\]
in \C, as $\nti$. Here $b$ and $\rho$ are the constants of Equation~\eqref{y-asym} and $l(t)$
denote the local time of a standard scaled Brownian excursion. 
\eth

In order to prove this result we have to show the following two theorems

\bth\label{th_fdd}
Let $b$, $\rho$, and $l_n(t)$ be as in Theorem~\ref{theo_1}, then for any $d$ and any choice of
fixed numbers $t_1,\dots, t_d$ the following limit theorem holds: 
$$
(l_n(t_1),\dots,l_n(t_d))\cw \frac{b\sqrt \rho}{2\sqrt2} \(l\left(\frac{b\sqrt \rho}{2\sqrt2}
\cdot t_1\right),\dots,l\left(\frac{b\sqrt \rho}{2\sqrt2}\cdot t_d\right)\),
$$
as $n\to\infty$. 
\eth

\bth\label{th_tight}
For all non-negative integers $n,r,h$ we have 
\begin{equation}
\E\left( L_n(r)-L_n(r+h)\right)^4 \le C\, h^2 n
\label{tight}
\end{equation}
where $C$ denotes some fixed positive constant. Consequently, the process
$l_n(t)$ is tight. 
\eth

\section{Combinatorial Setup}

In order to compute the distribution of the number of nodes in some given
levels in a tree of size $n$ we have to calculate 
the number $y_{k_1m_1k_2m_2\cdots k_dm_dn}$ of trees 
of size $n$ with $m_i$ nodes in level $k_i$, $i=1,\dots,d$ and normalize by
$y_n$. 

Therefore we introduce the generating functions $y_k(x,u)$ defined 
by the recurrence relation 
\begin{align}
y_0(x,u)&=uy(x) \nonumber \\
y_{k+1}(x,u)&=x\exp\(\sum_{i\ge 1}\frac{y_k(x^i,u^i)}i\),\quad k\ge 0.
\label{rec}
\end{align}

The function $y_k(x,u)$ represents trees where the nodes in level $k$ are marked (and counted by
$u$). If we want to look at two levels at once, say $k$ and $\ell$, then we have to take trees
with height at most $k$ and substitute the leaves in level $k$ by trees 
with all nodes at level $\ell-k$ marked (counted by $v$) and marking their roots as well (counted
by $u$). This leads to the generating function 
$y_{k,\ell}(x,u,v)=\tilde y_{k,\ell-k}(x,u,v)$ satisfying the recurrence relation 
\begin{align}
\tilde y_{0,\ell}(x,u,v)&=uy_{\ell}(x,v)  \nonumber \\
\tilde y_{k+1,\ell}(x,u,v)&=x\exp\(\sum_{i\ge 1}\frac{\tilde y_{k,\ell}(x^i,u^i,v^i)}i\),\quad
k\ge 0.
\label{tightrec}
\end{align}

In general we get therefore 
\begin{align*} 
y_{k_1,\dots,k_d}(x,u_1,\dots,u_d)&= \sum_{m_1,\dots,m_d,n\ge0}y_{k_1m_1k_2m_2\cdots k_dm_dn}
u_1^{m_1}\cdots u_d^{m_d}x^n \\
&=\tilde y_{k_1,k_2-k_1,k_3-k_2,\dots,
k_d-k_{d-1}}(x,u_1,\dots,u_d)
\end{align*} 
where
\begin{align*} 
\tilde y_{0,\ell_2,\dots,\ell_d}(x,u_1,\dots,u_d)&=u_1\tilde 
y_{\ell_2,\dots,\ell_d}(x,u_2,\dots,u_d) \\
\tilde y_{k+1,\ell_2,\dots,\ell_d}(x,u_1,\dots,u_d)&=x\exp\(\sum_{i\ge 1}\frac{\tilde
y_{k,\ell_2,\dots,\ell_d}(x^i,u_2^i,\dots,u_d^i)}i\), \quad k\ge 0.
\end{align*} 
The coefficients of these function are related to the process $L_n(t)$ (see \eqref{L_n} and the
lines before) by
$$
\p{L_n(k_1)=m_1,L_n(k_1+k_2)=m_2,\dots,L_n(k_1+k_2+\cdots+k_d)=m_d}=
\frac{y_{k_1m_1k_2m_2\cdots k_dm_dn}}{y_n}
$$
where the $k_i$ are integers and the probability space of the measure $\mathbf P$ is the set of
P\'olya trees with $n$ vertices equipped with the uniform distribution.

As claimed in Theorem~\ref{theo_1}, the process $l_n(t)=\rcp{\sqrt n}L_n(t)$ converges weakly to
Brownian excursion local time. From \cite{Hoog82} (cf. \cite{CoHo81,DG97} as well) we know that the
characteristic function $\phi(t)$ of the total local time of a standard Brownian excursion at
level $\kappa$ is 
\begin{equation} \label{local_time_int}
\phi(t)=1+\frac{\sqrt2}{\sqrt\pi}\int\limits_\gamma \frac{t\sqrt{-s}\exp(-\kappa\sqrt{-2s}\,)}
{\sqrt{-s}\exp(\kappa\sqrt{-2s}\,)-it\sqrt2 \sinh\(\kappa\sqrt{-2s}\,\)} e^{-s}\,ds
\end{equation} 
where $\gamma=(c-i\infty,c+i\infty)$ with some arbitrary $c<0$. The characteristic function of the 
joint distribution of the local time at several levels $\kappa_1,\dots,\kappa_d$ was computed in
\cite{DG97} (for $d=2$ already in \cite{CoHo81} albeit written down in a form which does not
exhibit the recursive structure) and is given by 

\begin{equation}\label{mvloctime}
\phi_{\kappa_1\dots\kappa_d}(t_1,\dots,t_d)= 1+\frac{\sqrt2}{i\sqrt{\pi}}
\int_\gamma f_{\kappa_1,\dots,\kappa_d}(x,t_1,\dots,t_d)e^{-x}\,dx,
\end{equation}
where
\begin{align*}
&f_{\kappa_1,\dots,\kappa_p}(x,t_1,\dots,t_d)=    
\Psi_{\kappa_1}(x,
t_1+\Psi_{\kappa_2-\kappa_1}(\dots\Psi_{\kappa_{d-1}-\kappa_{d-2}}(x,t_{d-1}+
\Psi_{\kappa_d-\kappa_{d-1}}(x,t_d))\cdots))
\end{align*}
with
$$
\Psi_{\kappa}(x,t)=\frac{it\sqrt{-x}\exp(-\kappa\sqrt{-2x}\,)
}{\sqrt{-x}\exp(\kappa\sqrt{-2x}\,)-it\sqrt2
\sinh\left(\kappa\sqrt{-2x}\right)}.
$$
For further studies of this distribution see \cite{GL99, Pi99, Ta91b}.

In order to 
show the weak limit theorem we have to show pointwise convergence of the characteristic function
$\phi_{k_1,\cdots, k_d,n}(t_1,\dots,t_d)$ of the joint distribution of $\rcp{\sqrt
n}L_n(k_1),\dots,\rcp{\sqrt n}L_n(k_d)$ to the corresponding characteristic function of the local
time in some interval containing zero. We have 
$$
\phi_{k_1,\cdots,k_d,n}(t_1,\dots,t_d)= 
\rcp{y_n}[x^n]y_{k_1,\dots,k_d}\left(x,e^{it_1/\sqrt n}, \dots
e^{it_d/\sqrt n}\right). 
$$

This coefficient will be calculated asymptotically by singularity analysis (see \cite{FO90}) for
$k_j=\lfloor \kappa_j\sqrt n\rfloor$. Thus knowing 
the local behaviour of $y_k(x,u)$ near its dominant singularity is the crucial step in
proving Theorem~\ref{theo_1}. This is provided by the following theorem will be the crucial step of the proof. 

\bth\label{y_k_lok}
Set $w_k(x,u)=y_k(x,u)-y(x)$. Let $x=\rho\left(1+\frac sn\right)$, $u=e^{it/\sqrt n}$, and
$k=\lfloor\kappa\sqrt n\rfloor$. 
Moreover, assume that $|\arg s|\ge \theta>0$ and, as $\nti$, we have
$s=\Ord{\log^2n}$ 
whereas $t$ and $\kappa$ are fixed. 
Then $w_k(x,u)$ admits the local representation 
\begin{align} 
w_k(x,u)
&\sim \frac{b^2\rho}{2\sqrt n}\cdot\frac{it\sqrt{-s}\exp\(-\kappa b\sqrt{-\rho
s}\,\)}{\sqrt{-s} -\frac{itb\sqrt\rho}4\(1-\exp\(-\kappa b \sqrt{-\rho
s}\,\)\)}
\label{locexp} \\
&=\frac{b\sqrt{2\rho}}{\sqrt n}\Psi_{\kappa b \sqrt\rho/(2\sqrt 2)}
\left( s, \frac{i t b\sqrt \rho}{2\sqrt 2}\right)
\end{align}
uniformly for $k=\Ord{\sqrt n}$.
\eth

\noindent
The proof is deferred to the next section.

Note that Theorem~\ref{y_k_lok} implies Theorem~\ref{th_fdd} for the case $d=1$. 
If we set $d=1$ in
Theorem~\ref{th_fdd}, then we have 
\begin{align*}
\phi_{k;n}(t)&=\rcp{y_n}[x^n]y_k\left(x,e^{it/\sqrt n}\right) \nonumber \\
&=\rcp{2\pi iy_n}\int_{\Gamma}y_k\left(x,e^{it/\sqrt n}\)\frac{dx}{x^{n+1}} \\
&=1+\rcp{2\pi iy_n} \int_{\Gamma}w_k\left(x,e^{it/\sqrt n}\)\frac{dx}{x^{n+1}}
\end{align*}
where $\Gamma$ is a suitable closed contour encircling the origin. Using Theorem~\ref{y_k_lok} it
is easy to show that $\phi_{\kappa\sqrt n;n}(t)$ converges to $\phi_{\kappa b \sqrt\rho/(2\sqrt
2)}(t b\sqrt \rho/2\sqrt 2)$ as desired. 

The higher dimensional case is more involved, but relies 
on the same principles. A complete proof is given in Section~\ref{fdds}.

\section{The Local Behaviour of $y_k(x,u)$ -- Proof of Theorem~\ref{y_k_lok}}\label{mainrec}

\subsection{Notation} We will provide some frequently used notations now. 

We will study the local behaviour of $y_k$ by analyzing the quantity  
$$
w_k(x,u)=y_k(x,u)-y(x)
$$
which frequently involves the term
$$
\Sigma_k(x,u):=\sum_{i\ge 2}\frac{w_k(x^i,u^i)}i.
$$
Furthermore, estimates of the partial derivatives 
$$
\gamma_k(x,u)=\pdiff{}{u}{}y_k(x,u) \mbox{ and } \gamma_k^{[i]}(x,u)=\pdiff{}{u}{i}y_k(x,u), \qquad
i\ge 2,
$$
will be needed.

The asymptotic analysis of $w_k$ (resp. $y_k$) enables us to apply Cauchy's integral formula and 
get the coefficients of $y_k(x,u)$ asymptotically (see the proof of Theorem~\ref{th_fdd} in the
next section) which eventually leads to an integral of the form
\eqref{local_time_int}. Therefore estimates for $y(x)$, provided in the next lemma, and the other
functions appearing in our analysis are needed. The estimates will be valid in various
domains. Therefore let us introduce  
\begin{align}\label{eqDelta}
\Delta & = \{ x \in \mathbb{C}: |x|< \rho + \eta,\
|\arg(x-\rho)| > \theta\}, \\
\Delta_\eps & =\{ x \in \mathbb{C}: |x-\rho|<\eps,\
|\arg(x-\rho)| > \theta \} \label{deltaeps} \\
\Theta & = \{ x \in \mathbb{C}: |x|< \rho + \eta,\
|\arg(x-\rho)| \neq 0\} \nonumber \\
\Xi_k & = \{ u \in \mathbb{C}: |u|\le 1, k |u-1|\le \tilde\eta \} \label{xik}
\end{align}
with $\eps,\eta,\tilde \eta>0$ and $0< \theta < \frac{\pi}2$. 

\brem 
In all the arguments which we will use in the
following proofs it is always assumed that $\eps$ and $\eta$ are sufficiently small even if it is
not explicitly mentioned. 
\erem

\begin{figure}
\begin{center}
$\begin{array}{ccc}
\epsfig{figure=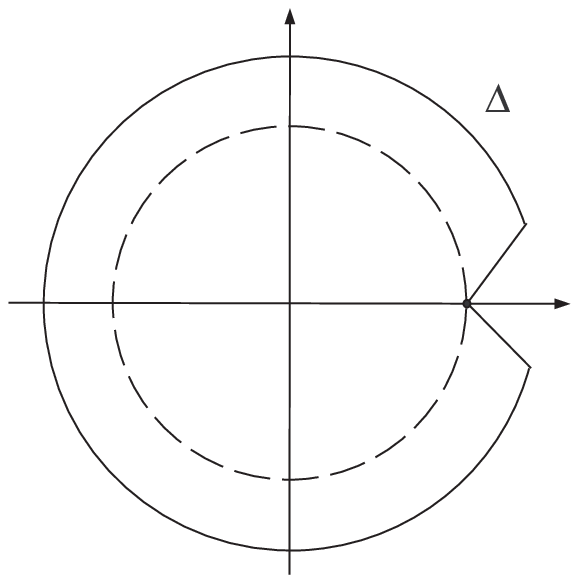,width=4cm}
&
\epsfig{figure=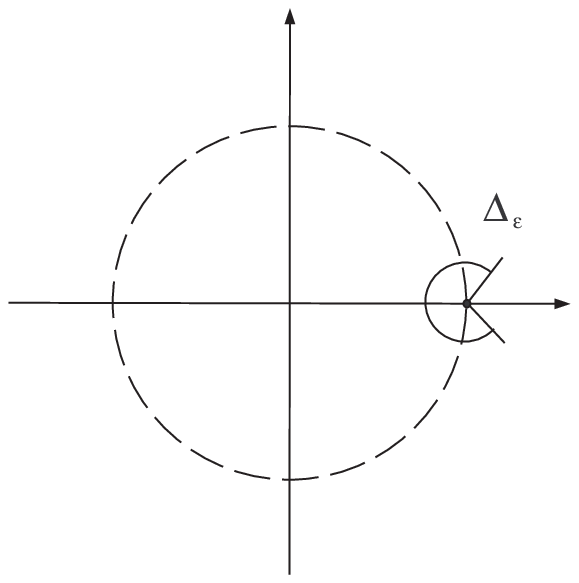,width=4cm}
&
\epsfig{figure=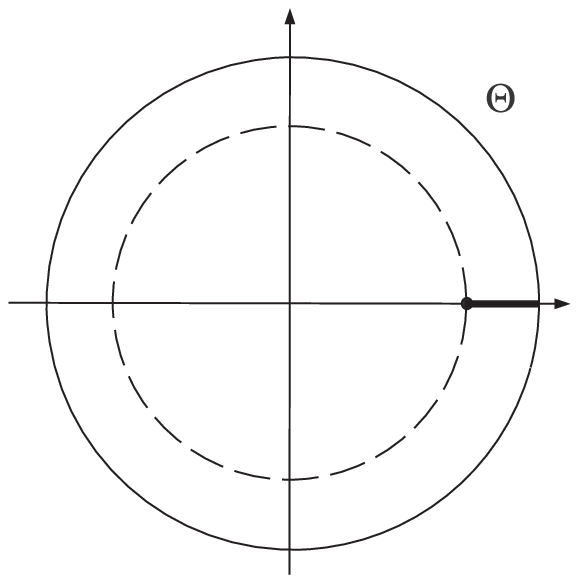,width=4cm}
\end{array}$
\end{center}
\caption{The domains $\Delta$, $\Delta_\eps$, and $\Theta$. For analyticity arguments we will need
the domain $\Delta$. Sometimes, the proof goes without change for the larger domain $\Theta$ as
well. In that case we will use $\Theta$ for the sake of more generality. For arguments where the
asymptotic behaviour near the singularity is important the domain $\Delta_\eps$ is used. The
domains $\Xi_k$ (not depicted here) are needed whenever the asymptotic behaviour for $u\sim 1$ is
considered. Here uniformity is often important such that a universal $\eps$-neighbourhood of $u=1$
which is independent of $k$ does not serve our needs.}\label{figure}
\end{figure}

\subsection{Analysis of the local behaviour of $y_k(x,u)$}

Obviously, $w_k(x,1)\equiv 0$. Since $y_k(x,u)$ represents the set of trees where the vertices of
level $k$ are marked, we expect that $\lim_{k\to\infty} y_k(x,u)=y(x)$ inside the domain of
convergence. This is not obvious, but follows from what we derive in the sequel. We start with a
useful property of $y(x)$.

\bl\label{lem_2}
Provided that $\eta$ in \eqref{eqDelta} is sufficiently small, 
the generating function $y(x)$ has the following properties:
\begin{enumerate}
\item[a)] For $x\in \Delta$ we have $|y(x)|\le 1$. Equality holds only for $x=\rho$.
\item[b)] Let $x=\rho\(1-\frac{1+it}n\)$ and $|t|\le C\log^2n$ for some fixed $C>0$. 
Then there is a $c>0$ such that 
$$
|y(x)|\le 1-c\sqrt{\frac{\max(1,|t|)}{n}}.
$$
\item[c)] For $|x|\le \rho$ we have $|y(x)|\le y(|x|)\le 1$. Moreover, near $x=0$ 
the asymptotic relation $y(x)\sim x$ holds. 
\item[d)] There exists an $\eps>0$ such that 
\begin{equation} \label{y_lb}
|y(x)|\ge \min\(\frac{\eps}2,\frac{|x|}2\)
\end{equation} 
for all $x\in \Theta$. 
\end{enumerate}
\el

\bpf
The first statement, when restricted to $|x|\le \rho$, 
follows from the facts that $y(x)$ has only positive coefficients (except
$y_0=0$), $y(\rho)=1$ and there are no periodicities. Extension to $\Delta$ is easily established
by using (\ref{y-asym}) and continuity arguments. 

The second statement is an immediate 
consequence of the singular expansion \eqref{y-asym} of $y(x)$. 

The first inequality of the third statement follows from the positivity of the coefficients $y_n$. 
The same fact also implies that $y(x)$ is strictly increasing in the interval $[0,\rho]$ and
therefore bounded by $y(\rho)=1$. The asymptotic relation near 0 follows from $y_0=0$ and $y_1=1$
or from the functional equation \eqref{funeq} 

Finally, for proving the last statement we split the circle $|x|\le \rho$ into the smaller circle
$|x|<\eps$ (let us call it $C$) and the annulus $A:=\{x\: : \: \eps\le|x|\le\rho\}$. The function
$y(x)$ is analytic in $C\cup A$ except at $x=\rho$, but still continuous there. Since the annulus
$A$ is compact, $|y(x)|$ attains a minimum there. By the functional equation \eqref{funeq} this
minimum must be positive, since both factor on the right-hand side of \eqref{funeq} are nonzero. 
Since $y(x)\sim x$, as $x\to 0$, a sufficiently small $\eps$ guarantees that $\min_{|x|=\eps}
|y(x)|>\eps/2$. Since this inequality holds for every smaller $\eps$ as well, we can choose an
$\eps$ such that 
$$
\min_{\eps\le |x|\le \rho} |y(x)|>\frac{\eps}2 \quad\mbox{ and }\quad |y(x)|\ge \frac{|x|}2 \mbox{ for }
|x|<\eps
$$
hold. By continuity this can be extended to \eqref{y_lb}. 
\epf

Next we derive an \emph{a priori} estimate of $w_k$ in a small domain. 

\bl \label{lem_1}
Let $|x|\le \rho^2+\eps$ for sufficiently small $\eps$ and $|u|\le 1$. Then there exist a constant
$L$ with $0<L<1$ and a positive constant $C$ such that 
$$
|w_k(x,u)|\le C|u-1|\cdot|x|\cdot L^k
$$
for all non-negative integers $k$.
\el

\bpf
We first note that by using the recurrence relation \eqref{rec} we obtain
\begin{align} 
w_{k+1}(x,u)&=y_{k+1}(x,u)-y(x) \nonumber\\
&=x\exp\(\sum_{i\ge 1}\rcp iy_k(x^i,u^i)\)-y(x) \nonumber \\
&=y(x)\(\exp\(w_k(x,u)+\sum_{i\ge 2}\frac{w_k(x^i,u^i)}i\)-1\) \label{wkrec}
\end{align} 
For $k=0$ we have $|w_0(x,u)|=|u-1|\cdot|y(x)|\le C|u-1||x|$ since $y(x)=\Ord{x}$ as $x\to 0$. We will then use the trivial inequality 
\begin{equation}\label{eqtrivinequ}
|e^x - 1| \le \frac{|x|}{1- \frac{|x|}2}
\end{equation}
for the induction steps. However, in order to apply this tool we
need some a-priori estimates.

Obviously we have for $|x|\le \rho$ and $|u|\le 1$
\[
|w_k(x,u)| \le 2 y(|x|)
\]
and consequently
\[
\left| \sum_{i\ge 1} \frac{w_k(x^i,u^i)}{i} \right| \le
2 \sum_{i\ge 1} \frac{y(|x|^i)}{i}  = 2 \log \frac{y(|x|)}{|x|}.
\]
Since the function $\frac{y(x)}{x}$ is convex for $0\le x\le \rho$
and $y(\rho) = 1$ we get $\frac{y(|x|)}{|x|}\le 1 + \frac{|x|}\rho$ because of the value of
$\rho$. Consequently
\[
\log \frac{y(|x|)}{|x|} \le \log \left( 1 + \frac{|x|}\rho \right)
\le \frac{|x|}\rho.
\]
Thus, if $|x|\le \rho^2 + \varepsilon$ (for a sufficiently small
$\varepsilon> 0$ we have
\[
\left| \sum_{i\ge 1} \frac{w_k(x^i,u^i)}{i} \right| \le 2 \rho
+ 2 \frac{\varepsilon}{\rho}.
\]
By using (\ref{eqtrivinequ}) we thus obtain
\[
\left| \exp\left( \sum_{i\ge 1} \frac{w_k(x^i,u^i)}{i} \right) - 1 \right|
\le \frac 1{ 1 - \rho - \frac {\varepsilon}{\rho} }
\left| \sum_{i\ge 1} \frac{w_k(x^i,u^i)}{i} \right|.
\]
Therefore, if we assume that we already know 
$|w_k(x,u)| \le C |u-1| |x| L^k$ (for $|x|\le \rho^2+\eps$, $|u|\le 1$
and some $L$ with $0< L < 1$) then we also get
\begin{align*}
|w_{k+1}(x,u)| &\le |y(x)| \cdot
 \left| \exp\left( \sum_{i\ge 1} \frac{w_k(x^i,u^i)}{i} \right) - 1 \right| \\
 &\le \frac {|y(x)|}{ 1 - \rho - \frac {\varepsilon}{\rho} }
 \left| \sum_{i\ge 1} \frac{w_k(x^i,u^i)}{i} \right| \\
 &\le \frac {|y(x)|}{ 1 - \rho - \frac {\varepsilon}{\rho} }
  C L^k\sum_{i\ge 1} \frac{|u^i-1|}i |x|^i \\
 & \le \frac {|y(x)|}{ 1 - \rho - \frac {\varepsilon}{\rho} }
 C L^k |u-1| \frac {|x|}{1-|x|}.
\end{align*}
By convexity we have $y(x) < x/\rho$ for $0< x < \rho$ and, thus,
there exists $\varepsilon> 0$ with
 $y(\rho^2+\varepsilon) \le \rho$.
Consequently we get for $|x|\le \rho^2 + \varepsilon$ the estimate
\[
|w_{k+1}(x,u)| \le C L' L^k |x| |u-1|
\]
holds where 
\[
L' = \frac {\rho}{\left( 1 - \rho - \frac {\varepsilon}{\rho}\right)
(1 - \rho^2 - \varepsilon) }. 
\]
The value of $L'$ is smaller than $1$ if $\varepsilon> 0$ is sufficiently small.
Thus, an induction proof works for $L = L'$.
\epf

\bcor\label{cor}
For $|u|\le 1$ and $|x|\le \rho+\eps$ ($\eps>0$ small enough) 
there is a positive constant $\tilde C$ such that (for all $k\ge
0$)
$$
|\Sigma_k(x,u)|\le \tilde C|u-1|L^k
$$
with the constant $L$ from the previous lemma.
\ecor

\bpf
We have 
\begin{align*} 
|\Sigma_k(x,u)|&\le \sum_{i\ge 2}\rcp i|w_k(x^i,u^i)|\le C \sum_{i\ge 2}\rcp i|u^i-1|\cdot|x|^iL^k
\\
&\le C|u-1|L^k\frac{|x|^2}{1-|x|}\le C|u-1|L^k\frac{1}{1-(\rho+\eps)}=\tilde C|u-1|L^k
\end{align*} 
\epf

\bcor\label{zcor}
Let $u\in\Xi_k$ and $x\in \Delta_\eps$. Then  
$$
\sum_{i\ge 2}\gamma_k(x^i,u^i)=\Ord{L^k}.
$$
\ecor

\bpf
The assertion follows immediately from the previous corollary and Taylor's theorem.
\epf

We will eventually need precise upper and lower bounds of $w_k(x,u)$ near its singularity. 
By Taylor's theorem $w_k(x,u)$ can be expressed in terms of $\gamma_k(x,u)$. Hence, let us
consider $\gamma_k(x,u)$ now.

\bl\label{Le7}
For $x\in\Theta$ (where $\eta>0$ is sufficiently small)
the functions $\gamma_k(x)$ can be represented as
\begin{equation} \label{gammakasym}
\gamma_k(x):=\gamma_k(x,1) = C_k(x)y(x)^k,
\end{equation} 
where the $C_k(x)$ form a sequence of analytic functions which converges uniformly 
to an analytic limit function $C(x)$ (for $x\in \Theta$) 
with convergence rate
\begin{equation} \label{Ckestim}
C_k(x) = C(x) + O(L^k),
\end{equation} 
for some $L$ with $0<L<1$. Furthermore we have $C(\rho) = \frac 12 b^2 \rho$.

There exist constants $\eps,\theta,\tilde \eta>0$ and $\theta<\frac\pi2$ such that 
\begin{equation}\label{eqgammaxuest}
|\gamma_k(x,u)| = O(|y(x)|^{k+1})
\end{equation}
uniformly for $x\in \Delta_\eps$ and $u\in \Xi_k$.
\el

\brem 
Note that the constants $\eps,\theta,\tilde \eta>0$ determine the domains $\Delta_\eps$ and
$\Xi_k$. Thus the theorem guarantees that there are suitable domains $\Delta_\eps$ and 
$\Xi_k$ of the shape \eqref{deltaeps} and \eqref{xik} where the estimates are valid.
\erem

\bpf
The first statement we have to show is that the functions $\gamma_k(x)$ are analytic functions in
$\Theta$. We prove this by induction. Obviously, the assertion holds for $k=0$ since
$\gamma_0(x)=y(x)$. Assume it is true for $k$. Then $\Gamma_k(x):=\sum_{i\ge 2} \gamma_k(x^i)$ is
analytic for $|x|<\sqrt\rho$ and hence also in $\Theta$. 
Using the recurrence relation of $y_k(x,u)$, Equation~\eqref{rec}, we get
\begin{align}
\gamma_{k+1}(x,u)&=\pdiff{}{u}{} x e^{y_k(x,u)+\Sigma_k(x,u)} \nonumber  \\
&=x e^{y_k(x,u)+\Sigma_k(x,u)} \sum_{i\ge 1}  \pdiff{}u{} y_k(x^i,u^i)u^{i-1}
\nonumber
\\
& = y_{k+1}(x,u) \sum_{i\ge 1} \gamma_k(x^i,u^i) u^{i-1}.
\label{eqgammaxurec}
\end{align}
This implies 
$$
\gamma_{k+1}(x)=y(x)\gamma_k(x)+y(x)\Gamma_k(x)
$$
which finally implies that $\gamma_{k+1}(x)$ is analytic in $\Theta$ as well. 

By solving this recurrence we obtain also the analyticity of
$C_k(x)=\frac{\gamma_k(x)}{y(x)^k}$ in $\Theta$. 
Furthermore, we will show that the sequence $(C_k(x))_{k\ge 0}$ has 
a uniform limit $C(x)$ which has the desired properties. 

Setting $u=1$ we can rewrite \eqref{eqgammaxurec} to 
\begin{equation} \label{ckestim}
C_{k+1}(x)y(x)^{k+2}=C_k(x)y(x)^{k+2}+ y(x) \(C_k(x^2)y(x^2)^{k+1}+ C_k(x^3)y(x^3)^{k+1}+\dots\).
\end{equation} 
resp.\ to
\begin{equation}\label{eqCkrec}
C_{k+1}(x) = \sum_{i\ge 1} C_k(x^i) \frac{y(x^i)^{k+1}}{y(x)^{k+1}}.
\end{equation}
Set 
\[
L_k := \sup_{x\in\Theta} 
\sum_{i\ge 2} \frac{|y(x^i)|^{k+1}}{|y(x)|^{k+1}}.
\]
If $\eta>0$ is sufficiently small then due to Lemma~\ref{lem_2} we have  
\[
\sup_{x\in\Theta} \frac{|y(x^i)|}{|y(x)|} < 1 \quad
\mbox{for all $i\ge 2$ and}\quad
\sup_{x\in\Theta}
 \frac{|y(x^i)|}{|y(x)|} = O(\overline L^i)
\]
for some $\overline L$ with $0< \overline L< 1$. Consequently we also get
\[
L_k = O(L^k)
\]
for some $L$ with $0< L< 1$ (actually we can choose $L=\overline L^2$). 
Thus, if we use the notation
$\displaystyle \| f\| := \sup_{x\in\Theta} |f(x)|$
then (\ref{eqCkrec}) yields 
\begin{equation}\label{eqCkrec2}
\|C_{k+1}\| \le \| C_k\| (1 + L_k)
\end{equation}
and also
\begin{equation}\label{eqCkrec3}
\| C_{k+1} - C_k\| \le \| C_k \| L_k.
\end{equation}
But (\ref{eqCkrec2}) implies that the functions
$C_k(x)$ are uniformly bounded in the given domain by
\[
\| C_k\| \le c_0 := \prod_{\ell\ge 1} (1+L_\ell). 
\]
Furthermore, (\ref{eqCkrec3}) guarantees the existence of a limit
$\lim_{k\to\infty} C_k(x) = C(x)$ which is analytic in $\Theta$; 
and we have uniform exponential convergence rate
\[
\| C_k - C\| \le c_0 \sum_{\ell \ge k} L_\ell = O(L^k).
\]
Hence, we get \eqref{gammakasym} as desired. 

Finally, note that (for $|x|\le \rho$)
\[
\sum_{k\ge 0} \gamma_k(x,1) = \sum_{n\ge 1} n y_n x^n=xy'(x).
\]
On the other hand, 
\[
\sum_{k\ge 0} \gamma_k(x,1) = \sum_{k\ge 0} (C(x)+\Ord{L^k})y(x)^{k+1}n= 
\frac{C(x)y(x)}{1-y(x)}+ \Ord{1}.
\] 
Since $C(x)$ is continuous in $\Theta$ this implies 
$$
C(\rho)=\lim_{x\to\rho} \frac{xy'(x)(1-y(x))}{y(x)}=\frac{b^2\rho}2
$$
where we used \eqref{y-asym}.

Let us turn to the second assertion. In order to obtain the upper bound (\ref{eqgammaxuest}) we
set for $\ell \le k$
\begin{equation} \label{Cl}
\overline C_\ell = \sup_{x\in\Delta_\eps,\ u\in \Xi_k} |\gamma_\ell(x,u) y(x)^{-\ell}|.
\end{equation} 
Observe that by Lemma~\ref{lem_2} $|\gamma_\ell(x,u)|\le \overline C_\ell \sup_{x\in\Delta_\eps}
|y(x)|^k\le \overline C_\ell$. Therefore, by Taylor's theorem and Corollary~\ref{cor} we
obtain 
\begin{align}
|y_{\ell+1}(x,u)| &\le |y(x)| \exp\left( \sum_{i\ge 1} \frac{|w_\ell(x^i,u^i)|}i \right)
\nonumber\\
& \le |y(x)| e^{\overline C_\ell|u-1| + O(L^\ell)}. \label{auxili}
\end{align}
By \eqref{eqgammaxurec} we have 
$$
\overline C_{\ell+1} =\sup_{x\in\Delta_\eps,\ u\in\Xi_k} \left|
\frac{y_{\ell+1}(x,u)}{y(x)^{\ell+1}} \right|\cdot \left|\sum_{i\ge 1}
\gamma_\ell(x^i,u^i)u^{i-1}\right|.
$$
Applying \eqref{auxili} and the estimate in Corollary~\ref{zcor} we get
\begin{align}
\overline C_{\ell+1} &\le \sup_{x\in\Delta_\eps,\ u\in\Xi_k} \left|\frac{e^{\overline
C_\ell|u-1|+\Ord{L^\ell}}}{y(x)^\ell}\right| \cdot\left|\gamma_\ell(x,u)+\sum_{i\ge
2}\gamma_\ell(x^i,u^i)u^{i-1}\right| \nonumber \\
&\le e^{\overline C_\ell\tilde\eta /k} \overline C_\ell \(1+\Ord{L^\ell}\) \label{eqDkrec} 
\end{align}

Set
\[
c_0 = \prod_{j\ge 0} (1+ O(L^j))
\]
and choose $\tilde\eta>0$ such that $e^{2c_0\tilde\eta }\le 2$.
We also choose $\eps\le \tilde\eta$ and $0< \theta < \frac \pi 2$
such that $|y(x)|\le 1$ for $|x-\rho|< \eps$ and $|\arg(x-\rho)|\ge \theta$.
If $k> 0$ is fixed, then it 
follows by induction that, provided that $|u-1|\le \tilde\eta /k$,
\[
\overline C_\ell \le \prod_{j< \ell} (1+O(L^j)) \cdot e^{2c_0\tilde\eta  \ell/k} \le 2 c_0 \qquad
(\ell \le k).
\]
This completes the proof of the lemma, since $\bar C_0=1$. 
\epf

The representation (\ref{gammakasym}) from Lemma~\ref{Le7}
gives us a first indication of the behaviour of $w_k(x,u)$
for $u$ close to $1$. We expect that
\begin{equation}\label{eqfirstappr}
w_k(x,u) \approx (u-1)\gamma_k(x,1) \sim (u-1)C(x) y(x)^k.
\end{equation}
This actually holds (up to constants) in a proper range for $u$ and $x$, 
although it is only partially true 
in the range of interest (see Theorem~\ref{y_k_lok}). 

In order to make this more precise we derive estimates for
the second derivatives $\gamma_k^{[2]}(x,u)$.

\bl\label{Le7-2}
Suppose that $|x|\le \rho-\eta$ for some $\eta > 0$ and $|u|\le 1$. Then 
uniformly
\begin{equation}\label{eqgammaxu2est1}
\gamma_k^{[2]}(x,u) = O(y(|x|)^{k+1}).
\end{equation}
There also exist constants $\eps,\theta,\tilde \eta>0$ such that
\begin{equation}\label{eqgammaxu2est}
\gamma_k^{[2]}(x,u) = O( k\, |y(x)|^{k+1})
\end{equation}
uniformly for $u\in\Xi_k$ and $x\in\Delta_\eps$.
\el

\bpf
By definition we have the recurrence (compare with \eqref{eqgammaxurec}) 
\begin{align}
\gamma_{k+1}^{[2]}(x,u) &= y_{k+1}(x,u) \sum_{i\ge 1} i\gamma_k^{[2]}(x^i,u^i)u^{2i-2} 
 \nonumber \\
&+ y_{k+1}(x,u) \left(\sum_{i\ge 1} \gamma_k(x^i,u^i) u^{i-1} \right)^2 
 \label{eqgammak2rec-0}\\
&+ 
y_{k+1}(x,u) \sum_{i\ge 2} (i-1)\gamma_k(x^i,u^i) u^{i-2} \nonumber
\end{align}
with initial condition $\gamma_{0}^{[2]}(x) = 0$.

First suppose that $|x|\le \rho-\eta$ for some $\eta> 0$ and $|u|\le 1$.
Then we have $|\gamma_{k}^{[2]}(x,u)| \le \gamma_{k}^{[2]}(|x|,1)$.
Thus, in this case it is sufficient to consider non-negative real $x\le \rho-\eta$. 
We proceed by induction. 
Suppose that we already know that $\gamma_{k}^{[2]}(x):=\gamma_{k}^{[2]}(x,1)
 \le D_k y(x)^{k+1}$ (where $D_0 = 0$).
Then we get from (\ref{eqgammak2rec-0}) and the already known bound
$\gamma_k(x,1) \le C y(x)^k$ from Lemma~\ref{Le7} the upper bound
\begin{align*}
\gamma_{k+1}^{[2]}(x) &\le  D_k y(x)^{k+2} + 
D_k y(x)^{k+2} \sum_{i\ge 2} i \frac{y(x^i)^{k+1}}{y(x)^{k+1}}\\
&\qquad+ y(x)\(C^2 y(x)^{2k}  + \sum_{i\ge 2} \gamma_k(x^i,1)\)^2 \\
&\qquad+ C  y(x)^{k+2} \sum_{i\ge 2}(i-1)  \frac{y(x^i)^{k+1}}{y(x)^{k+1}} \\
&\le y(x)^{k+2}\left( D_k( 1 + O(L^k))  + C^2 y(\rho-\eta)^k + O(L^k) \right)
\end{align*}
where we used Corollary~\ref{zcor} in the last step. 
Consequently we can set
\[
D_{k+1} = D_k( 1 + O(L^k))  + C^2 y(\rho-\eta)^k + O(L^k)
\]
and obtain that $D_k = \Ord{1}$ as $k\to\infty$ which proves
(\ref{eqgammaxu2est1}).

Next choose $\eps,\theta,\tilde \eta$ as in the proof of Lemma~\ref{Le7} and 
set (for $\ell \le k$)
\[
\overline D_\ell = \sup_{x\in\Delta_\eps,\ u\in\Xi_k} 
\left|\gamma_{\ell}^{[2]}(x,u) y(x)^{-\ell-1}\right|.
\]
By the same reasoning as in 
the proof of Lemma~\ref{Le7}, where we use the already proved
bound $|\gamma_\ell(x,u)\le \overline C |y(x)|^{\ell+1}$, we obtain 
\[
\overline D_{\ell+1} \le \overline D_\ell\,e^{\tilde\eta C/k} 
(1 + O(L^\ell)) 
+ C^2 \,e^{\tilde\eta C/k}  + O(L^\ell),
\]
that is, we have
\[
\overline D_{\ell+1} \le \alpha_\ell \overline D_\ell + \beta_\ell
\]
with $\alpha_\ell = e^{\tilde\eta C/k} 
(1 + O(L^\ell))$ and $\beta_\ell = C^2 \,e^{\tilde\eta C/k}  + O(L^\ell)$. 
Hence we get
\begin{align*}
\overline D_k &\le \sum_{j=0}^{k-1} \beta_j \prod_{i=j+1}^{k-1} \alpha_i\\
&\le  k \, \max_j \beta_j \, e^{\tilde\eta C} \, \prod_{\ell\ge 0} (1 + O(L^\ell)) \\
& = O(k).
\end{align*}
This completes the proof of (\ref{eqgammaxu2est}).
\epf

Using the estimates for $\gamma_k(x,u)$ and $\gamma_k^{[2]}(x,u)$
we derive the following representations for $w_k(x,u)$ and $\Sigma_k(x,u)$.

\bl\label{Le7-3}
Let $\eps,\theta,\tilde\eta$ and $C_k(x) = \gamma_k(x,1)/y(x)^{k+1}$ as in Lemma~\ref{Le7}. Then
we have 
\begin{equation}\label{eqwkfirstrep}
w_k(x,u) = C_k(x)(u-1) y(x)^{k+1}\left( 1 + O(k|u-1|)\right),
\end{equation}
uniformly for $u\in \Xi_k$ and $x\in\Delta_\eps$.

Furthermore we have for $|x|\le \rho+\eta$ (for some $\eta>0$) 
and $|u|\le 1$
\begin{equation}\label{eqSigmakfirstrep}
\Sigma_k(x,u) = \tilde C_k(x)(u-1) y(x^2)^{k+1}+ O\left(|u-1|^2 y(|x|^2)^k\right),
\end{equation}
where the analytic functions $\tilde C_k(x)$ are given by
\begin{equation}\label{eqtildeCkrep}
\tilde C_k(x) = \sum_{i\ge 2} C_k(x^i) 
\left( \frac{y(x^i)}{y(x^2)} \right)^{k+1}.
\end{equation}
They have a uniform limit
$\tilde C(x)$ with convergence rate
\[
\tilde C_k(x) = \tilde C(x) + O(L^k)
\]
for some constant $L$ with $0< L < 1$.
\el

\bpf
The first relation (\ref{eqwkfirstrep}) follows from 
Lemma~\ref{Le7}, Lemma~\ref{Le7-2} and Taylor's theorem. 

In order to prove (\ref{eqSigmakfirstrep}) we first note that
$|x^i|\le \rho-\eta$ for $i\ge 2$ and $|x|\le \rho+\eta$
(if $\eta>0$ is sufficiently small). Hence, by a second use
of Taylor's theorem we get uniformly 
\[
w_k(x^i,u^i) = C_k(x^i)(u^i-1)y(x^i)^{k+1}+ O\left(|u^i-1|^2 y(|x^i|)^{k+1}\right)
\]
and consequently
\begin{align*}
\Sigma_k(x,u) &= \sum_{i\ge 2} \frac 1i\, C_k(x^i)(u^i-1)y(x^i)^{k+1} + 
O\left( |u-1|^2 y(|x|^2)^k \right) \\
&= (u-1) \tilde C_k(x) y(x^2)^{k+1} +O\left( |u-1|^2 y(|x|^2)^k \right).
\end{align*}
Here we have used the property that the sum
\[
\sum_{i\ge 2}  C_k(x^i)\frac{u^i-1}{i(u-1)} \frac{y(x^i)^{k+1}}{y(x^2)^{k+1}}
\]
represents an analytic function in $x$ and $u$, since due to $|x|<\rho+\eps<1$ we have
$x^i\to 0$ and therefore \eqref{Reihe} implies $y(x^i)\sim x^i$. Finally, since
$C_k(x) = C(x)+ O(L^k)$, it also follows that $\tilde C_k(x)$ has a limit
$\tilde C(x)$ and the same order of convergence.
\epf

With these auxiliary results we are able to get a precise result for $w_k(x,u)$.	
\bl\label{lem_4}
For $u\in \Xi_k$ and $x\in\Delta_\eps$ ($\eps,\theta,\tilde\eta$ as in Lemma~\ref{Le7}) we have  
\begin{equation} \label{wkfirst}
 w_k(x,u) =\frac{C_k(x)w_0(x,u)\, y(x)^{k}}{1-
\frac 12 C_k(x)w_0(x,u) \frac{1-y(x)^k}{1-y(x)}+\Ord{|u-1|}
}. 
\nonumber
\end{equation} 
\el

\bpf
Since by Lemma~\ref{Le7-3}  
$\Sigma_k(x,u) = O(w_k(x,u))$, we 
observe that $w_{k}(x,u)$ satisfies the recurrence relation 
(we omit the arguments now)
\begin{align*}
w_{k+1} &= y\left( e^{w_k+ \Sigma_k} -1 \right)\\
&= y\left( w_k + \frac{w_k^2}2 + \Sigma_k + O(w_k^3) + O(\Sigma_k^2)\right)\\
&=yw_k\(1+\frac{w_k}2+\Ord{w_k^2}+ \Ord{\Sigma_k}\)
\left( 1+ \frac {\Sigma_k}{w_k} \right).
\end{align*}
Equivalently, we have
\begin{align*}
\frac{y}{w_{k+1}} + \frac{y\,\Sigma_k}{w_kw_{k+1}} &=
\frac 1{w_k}\left( 1 - \frac{w_k}2+\Ord{w_k^2}+ \Ord{\Sigma_k}\right)\\
&=\frac 1{w_k} - \frac 12 + \Ord{w_k} + \Ord{\frac{\Sigma_k}{w_k}},
\end{align*}
and consequently
\[
\frac{y^{k+1}}{w_{k+1}} = 
\frac {y^k}{w_k} - \frac{\Sigma_k y^{k+1}}{w_kw_{k+1}}
- \frac 12 y^k + \Ord{w_k y^k} + \Ord{\frac{\Sigma_k y^k}{w_k}}.
\]
Thus we get by recurrence
$$
\frac{y^k}{w_k} =\rcp{w_0}-\sum_{\ell=0}^{k-1}\frac{\Sigma_\ell}
{w_\ell w_{\ell+1}}y^{\ell+1}
-\rcp2\frac{1-y^k}{1-y} +
\Ord{\frac{1- L^k}{1- L}}+\Ord{w_0\frac{1-y^{2k}}{1-y^2}}.
$$
Now we use again Lemma~\ref{Le7-3} to obtain
\begin{align*}
w_0\sum_{\ell=0}^{k-1}\frac{\Sigma_\ell}{w_\ell w_{\ell+1}}y^{\ell +1}
&= \sum_{\ell=0}^{k-1} \frac{\tilde C_\ell(x) y(x^2)^{\ell+1} + 
O\left(|u-1| y(|x|^2)^\ell\right) }
{C_\ell(x) C_{\ell+1}(x) y(x)^{\ell+1} (1 + O(\ell|u-1|)) }\\
&=\sum_{\ell = 0}^{k-1} \frac{\tilde C_\ell(x)}{C_\ell(x) C_{\ell+1}(x)}
\frac {y(x^2)^{\ell+1}}{y(x)^{\ell+1}}  + O(|u-1|) \\
&= c_k(x) + O(u-1)
\end{align*}
where $c_k(x)$ denotes the sum in the penultimate line above.
Observe, too, that $w_0\frac{1-y^{2k}}{1-y^2} = O(1)$, if $k|u-1| \le \tilde\eta $.
Hence we obtain the representation
\begin{equation}\label{eqcompare}
w_k = \frac{w_0y^k}{1-c_k(x) - \frac{w_0}2 \frac{1-y^k}{1-y}+ \Ord{|u-1|}}.
\end{equation}
Thus, it remains to verify that $1-c_k(x) = 1/C_k(x)$.
By using (\ref{eqCkrec}) and (\ref{eqtildeCkrep}) it follows that 
\begin{equation}\label{eqtildeCkCkrel}
\tilde C_k(x) = \sum_{i\ge 2} C_k(x^i) 
\left( \frac{y(x^i)}{y(x^2)} \right)^{k+1} = 
(C_{k+1}(x)- C_k(x)) \left( \frac{y(x)}{y(x^2)} \right)^{k+1}
\end{equation}
and consequently by telescoping
\[
c_k(x) = \sum_{\ell=0}^{k-1} \frac{C_{\ell+1}(x) - C_\ell(x)}
{C_\ell(x) C_{\ell+1}(x)} = \frac 1{C_0(x)} - \frac 1{C_k(x)} 
= 1 -  \frac 1{C_k(x)}.
\]
Alternatively we can compare (\ref{eqcompare}) with (\ref{eqwkfirstrep})
for $u\to 1$ which also shows $1-c_k(x) = 1/C_k(x)$.
This completes the proof of the lemma.
\epf

The proof of Theorem~\ref{y_k_lok} is now immediate.
We substitute $x=\rho\left(1+\frac sn\right)$ (where $s/n\to 0$), 
$u=e^{it/\sqrt n}$ and set $k = \kappa \sqrt n$.
We also use the local expansion
$y(x) = 1 - b\sqrt \rho \sqrt{1 - x/\rho} + O(|x-\rho|)$. That leads to
\[
y(x)^k = \exp\left( - \kappa b\sqrt{-\rho s} \right) \left(
1 + O\left( \frac \kappa{\sqrt n} \right) \right).
\]
Finally, since the functions $C_k(x)$ are continuous and uniformly
convergent to $C(x)$, they are also uniformly continuous and, thus,
\begin{equation} \label{Ckapprox} 
C_k(x) \sim C(\rho) = \frac 12 b^2 \rho. 
\end{equation} 
Altogether this
leads to 
\begin{align*}
&\frac{w_0(x,u)\, y(x)^{k}}{\frac 1{C_k(x)}-
\frac{w_0(x,u)}2 \frac{1-y(x)^k}{1-y(x)}+\Ord{|u-1|} }\\
&\quad=
\frac{C_k(x)(u-1)(1-y(x))\, y(x)^{k+1}}{1-y(x)-
C_k(x)\frac{w_0(x,u)}2 (1-y(x)^k)+\Ord{|u-1|\cdot|1-y(x)|} }\\
&\quad \sim \frac{b^2\rho}{2\sqrt n}\cdot\frac{it \sqrt{-s}
\exp\(-\frac 12\kappa b\sqrt{-\rho s}\,\)}
{\sqrt{-s}\exp\(\frac 12\kappa b\sqrt{-\rho s}\,\) -\frac{it b\sqrt\rho}2
\sinh\(\frac 12\kappa b\sqrt{-\rho s}\,\)}
\end{align*}
as proposed.

\section{The Finite Dimensional Limiting Distributions -- Proof of
Theorem~\ref{th_fdd}}\label{fdds}

For $d=1$ (in Theorem~\ref{th_fdd}) we have 
\begin{align}
\phi_{k,n}(t)&=\rcp{y_n}[x^n]y_k\left(x,e^{it/\sqrt n}\right) \nonumber \\
&=\rcp{2\pi iy_n}\int_{\Gamma}y_k\left(x,e^{it/\sqrt n}\)\frac{dx}{x^{n+1}} \label{cif}
\end{align}
where the contour $\Gamma=\gamma\cup\Gamma'$ consists of a line 
$$
\gamma=\{x=\rho\left(1-\frac{\sigma+i\tau}n\right)\,\mid\, -C\log^2n\le \tau \le C\log^2n\}
$$
with an arbitrarily chosen fixed constants $C>0$ and $\sigma>0$, and $\Gamma'$ 
is a circular arc centered at the
origin and making $\Gamma$ a closed curve.  

The contribution of $\Gamma'$ is exponentially small since for $x\in\Gamma'$ we have
$\rcp{y_n}|x^{-n-1}|=\Ord{n^{3/2} e^{-\log^2n}}$ 
whereas $\left|y_k\left(x,e^{it/\sqrt n}\)\right|$ is bounded. 

If $x\in\gamma$, then the local expansion \eqref{locexp} is valid. 
Insertion into \eqref{cif}, using \eqref{yn}, and taking the limit for $\nti$ yields the 
characteristic function of
the distribution of $\frac{b\sqrt\rho}{2\sqrt2} l\(\frac{b\sqrt{\rho}}{2\sqrt2}\kappa\)$ as
desired.

Now we proceed with $d=2$. The computation of the two dimensional limiting distributions shows the
general lines of the proof. An iterative use of the techniques will eventually prove
Theorem~\ref{th_fdd}. We confine ourselves with the presentation of the case $d=2$.

We have to show 
\begin{equation} \label{bivar}
\left(\frac 1{\sqrt n} L_n (\kappa_1 \sqrt n), \,
\frac 1{\sqrt n} L_n (\kappa_2 \sqrt n)\right)
 \cw \left( \frac {b\sqrt \rho}{2\sqrt 2}\cdot
l\left( \frac {b\sqrt \rho}{2\sqrt 2} \kappa_1\right),\,
\frac {b\sqrt \rho}{2\sqrt 2}\cdot
l\left( \frac {b\sqrt \rho}{2\sqrt 2} \kappa_2\right)\right).
\end{equation} 
Since the characteristic function of the two dimensional distribution satisfies 
\begin{align} 
\phi_{k,k+h,n}(t_1,t_2)&=\rcp{y_n}[x^n]\tilde y_{k,h}\left(x,e^{it_1/\sqrt
n},e^{it_2/\sqrt n}\) \nonumber \\
&= \rcp{2\pi iy_n}\int_{\Gamma}\tilde y_{k,h}\left(x,e^{it_1/\sqrt
n},e^{it_2/\sqrt n}\)\frac{dx}{x^{n+1}} \nonumber \\
&=1+\rcp{2\pi iy_n}\int_{\Gamma} w_{k,k+h}\left(x,e^{it_1/\sqrt
n},e^{it_2/\sqrt n}\)\frac{dx}{x^{n+1}}, \label{cauchy}
\end{align} 
where 
$$
w_{k,k+h}(x,u_1,u_2)=y_{k,k+h}(x,u_1,u_2)-y(x)=\tilde y_{k,h}(x,u_1,u_2)-y(x), 
$$
we need to analyze the asymptotic behaviour of $w_{k,k+h}(x,u_1,u_2)$ 
for $k$ and $h$ proportional to $\sqrt n$. Furthermore, note that $y(x)$ and $\tilde
y_{k,h}(x,u_1,u_1)$ are analytic functions for $x\in \Delta$ and thus $w_{k,k+h}$ is bounded for
$x\in\Delta$. Hence the contribution of $\Gamma'$ to the Cauchy integral \eqref{cauchy} is
$\Ord{1/y_n \cdot e^{-c\log^2n}}$ with some suitable constant $c>0$ and therefore negligibly small.
Extending $\gamma$ to infinity, as it is required for the two-dimensional version of
\eqref{local_time_int}, again introduces an exponentially small and therefore negligible error. 

Thus it is sufficient to know the behaviour of $w_{k,k+h}\left(x,e^{it_1/\sqrt
n},e^{it_2/\sqrt n}\)$ for $x\in\gamma$.

Then, equation \eqref{bivar} follows from the following proposition and \eqref{yn}.

\bpr\label{fddprop}
Let $\kappa_2 > \kappa_1 > 0$ and $t_1,t_2$ be given with
$|\kappa_2 t_1|\le c$ and $|\kappa_2 t_2|\le c$. Furthermore, set $u_i = e^{it_i/\sqrt n}$ and
$k_i=\lfloor \kappa_i \sqrt n\rfloor$ for $i=1,2$. Define $s$ by $x=\rho\(1+\frac sn\)$.
Then we have
\begin{equation} \label{biv}
w_{k_1,k_2}(x,u_1,u_2) \sim \frac{b\sqrt{2\rho}}{\sqrt n}\,
\Psi_{\frac{\kappa_1 b\sqrt\rho}{2\sqrt 2}}
 \left( s, \frac{it_1 b \sqrt\rho}{2\sqrt 2} +
\Psi_{\frac{(\kappa_2-\kappa_1) b\sqrt\rho}{2\sqrt 2}}
 \left( s, \frac{it_2 b \sqrt\rho}{2\sqrt 2}  \right)
 \right), 
\end{equation} 
as $\nti$, uniformly for $x$, $u_1$ and $u_2$ such that $k_1|u_1-1|\le c$, $k_2|u_2-1|\le c$, and 
$s=\Ord{\log^2n}$ such that $\Re(s) = -\sigma$ where $\sigma>0$ and sufficiently large but fixed.

\epr

\bpf
Note that $w_{k_1,k_2}(x,1,1) = 0$, $w_{k_1,k_2}(x,u_1,1) = w_{k_1}(x,u_1)$, and 
$w_{k_1,k_2}(x,1,u_2) = w_{k_2}(x,u_2)$. Thus, the first
derivatives are given by
\begin{align*}
\left[\frac{\partial}{\partial u_1} w_{k_1,k_2}(x,u_1,u_2)
\right]_{u_1=u_2 = 1} &= \gamma_{k_1}(x,1),\\
\left[\frac{\partial}{\partial u_2} w_{k_1,k_2}(x,u_1,u_2)
\right]_{u_1=u_2 = 1} &= \gamma_{k_2}(x,1).
\end{align*}
It is also possible to get bounds for the second derivatives 
of the form $O(k_2 y(x)^{k_1})$, if $x$, $u_1$, and $u_2$ are in the domain given in the assertion
above. Hence, we can approximate $w_{k_1,k_2}(x,u_1,u_2)$ by
\begin{align}
w_{k_1,k_2}(x,u_1,u_2) &= C_{k_1}(x)(u_1-1) y(x)^{k_1+1} +
 C_{k_2}(x)(u_2-1) y(x)^{k_2+1} \label{eqwk1k2rep}\\
&+ O\left( k_2 y(x)^{k_1} (|u_1-1|^2 +  |u_2-1|^2) \right). \nonumber
\end{align}
Similarly (and even more easily, compare with the proof of Lemma~\ref{Le7-3}) 
we obtain a representation for
\begin{align}
\Sigma_{k_1,k_2}(x,u_1,u_2) &= \sum_{i\ge 2}
\frac{w_{k_1,k_2}(x^i,u_1^i,u_2^i)} i \nonumber  \\
&= \tilde C_{k_1}(x)(u_1-1) y(x^2)^{k_1+1} + 
\tilde C_{k_2}(x)(u_2-1) y(x^2)^{k_2+1}  \nonumber  \\
& + O\left( y(|x|^2)^{k_1} |u_1-1|^2 +  y(|x|^2)^{k_2}|u_2-1|^2) \right).
\label{eqSigmak1k2rep} 
\end{align}
In order to identify the asymptotic main term of $\Sigma_{k_1,k_2}(x,u_1,u_2)$ note that 
$k_2 - k_1 \sim (\kappa_2-\kappa_1)\sqrt n$ and $|y(x^2)|<y(|x|^2)<1$ for $|x|\le \rho+\eta<1$. 
Moreover, observe that the terms $u_1-1$ and $u_2-1$ are proportional and hence 
\begin{equation} \label{auxil}
\frac{C_{k_2}(x)(u_2-1)y(x)^{k_2}}{C_{k_1}(x)(u_1-1)y(x)^{k_1}}\sim \frac{t_2}{t_1} \exp\left( - (
\kappa_2-\kappa_1 )b\sqrt{-\rho s} \right). 
\end{equation} 
The $\kappa_i$ and the $t_i$ are fixed, so we can choose $\sigma$ such that the right-hand side in
\eqref{auxil} is different from 1. This guarantees that the asymptotic main terms of
\eqref{eqwk1k2rep}, $C_{k_1}(x)(u_1-1)y(x)^{k_1}$ and $C_{k_2}(x)(u_2-1)y(x)^{k_2}$, 
do not cancel each other.

By using the same reasoning as in the proof of Lemma~\ref{lem_4}
we get the representation
\[
w_{k_1,k_2} = \frac{w_{0,k_2-k_1}y^{k_1}}
{1-f_{k_1} - \frac{w_{0,k_2-k_1}}2 \frac{1-y^{k_1}}{1-y} +
O(|u_1-1|+|u_2-1|) },
\]
where
\begin{align*}
f_{k_1} &= f_{k_1}(x,u_1,u_2) \\
&= w_{0,k_2-k_1}(x,u_1,u_2) \sum_{\ell = 0}^{k_1-1} 
\frac{\Sigma_{\ell,k_2-k_1+\ell}(x,u_1,u_2)\,  y(x)^{\ell+1}}
{w_{\ell,k_2-k_1+\ell}(x,u_1,u_2)w_{\ell+1,k_2-k_1+\ell+1}(x,u_1,u_2)}.
\end{align*}
Note that 
\begin{align*}
w_{0,k_2-k_1} &= u_1 t_{k_2-k_1}(x,u_2) - y(x) \\
&= (u_1-1)y(x) + u_1 w_{k_2-k_1}(x,u_2)\\
&= U + W,
\end{align*}
where $U$ and $W$ abbreviate $U = (u_1-1)y(x)$ and $W = u_1w_{k_2-k_1}(x,u_2)$.
The assumptions on $u_i$ and $k_i$ given in the statement of Proposition~\ref{fddprop} imply that
by Lemma~\ref{Le7-3} we have ($A\asymp B$ means that $A$ and $B$ have same order of magnitude)
$$
W \asymp u_1C_{k_2-k_1}(x)(u_2-1) y(x)^{k_2-k_1+1}=o(u_2-1)
$$
whereas $U\asymp u_1-1 \asymp u_2-1$. Thus we may safely assume that $|W|< \frac 12 |U|$, so 
that there is no cancellation.

Next by (\ref{eqwk1k2rep}) we have 
\begin{align*} 
w_{\ell,k_2-k_1+\ell} =&
C_\ell(x)(u_1-1)y(x)^{\ell+1}+C_{k_2-k_1+\ell}(x)(u_2-1)y(x)^{k_2-k_1\ell+1} \\
&+\Ord{(k_2-k_1\ell)y(x)^\ell(|u_1-1|^2+|u_2-1|^2)} \\
=&C_\ell (x)y(x)^\ell U+C_{k_2-k_1+\ell}(x)(u_2-1)y(x)^{k_2-k_1\ell+1}\\
&+\Ord{(k_2-k_1\ell)y(x)^\ell(|u_1-1|^2+|u_2-1|^2)}
\end{align*} 
But since $C_{k_2-k_1+\ell}(x)\sim C_{k_2-k_1}(x)$ by formula \eqref{Ckestim} and Lemma~\ref{Le7-3} relates the second term to $W$, we obtain 
\[
w_{\ell,k_2-k_1+\ell} = 
\left( C_\ell(x) U + W \right) y(x)^{\ell} 
\left( 1 (k_2-k_1)|u_2-1|) \right).
\]
Hence, $f_k$ can be approximated by 
(for simplicity we omit the error terms)
\begin{align*}
f_{k_1} &\sim U(U+W) 
\sum_{\ell= 0}^{k_1-1} \frac{\tilde C_\ell(x)\, (y(x^2)/y(x))^{\ell+1}}
{(C_\ell(x) U + W)(C_{\ell+1}(x) U + W)} \\
&= (U+W) \sum_{\ell= 0}^{k_1-1} \frac{(C_{\ell+1}(x) U + W)-(C_\ell(x) U + W) }
{(C_\ell(x) U + W)(C_{\ell+1}(x) U + W)} \\
&= (U+W) \left( \frac 1{U+W} - \frac 1{C_kU+W} \right)\\
&= 1 - \frac{U+W}{C_kU+W},
\end{align*}
where we have used the formula (\ref{eqtildeCkCkrel}) and telescoping.
Consequently, it follows that
\begin{align}
w_{k_1,k_2} &\sim \frac{(U+W)y^{k_1}}
{\frac{U+W}{C_kU+W} - \frac{U+W}2 \frac{1-y^{k_1}}{1-y} }\nonumber\\
&= \frac{(C_kU+W)y^{k_1}}
{1 - \frac{C_kU+W}2 \frac{1-y^{k_1}}{1-y} }.\label{approximendum}
\end{align}
Now we will approximate all the terms in \eqref{approximendum}. First recall $C_k(x) \sim
\frac{b^2\rho}2$ (formula \eqref{Ckapprox}). The assertion on $x=\rho\(1+\frac sn\)$ implies
$x\to\rho$ and hence $y(x)\to 1$. Thus we obtain 
$$
U = y(u_1-1) \sim \frac{it_1}{\sqrt n}.
$$
The asymptotic expansion \eqref{y-asym} and Theorem~\ref{y_k_lok} imply 
\begin{align*}
y^{k_1} &\sim \exp\left( - \frac 12 \kappa_1 b \sqrt{-\rho s} \right),\\
1-t &\sim b\sqrt \rho \sqrt{\frac sn}, \\
W &= u_1w_{k_2-k_1}(x,u_2) \sim \frac{b\sqrt{2\rho}}{\sqrt n}\,
\Psi_{\frac{(\kappa_2-\kappa_1) b\sqrt\rho}{2\sqrt 2}}
\left( s, \frac{it_2 b \sqrt\rho}{2\sqrt 2} \right) , 
\end{align*}
Applying all these approximations we finally get \eqref{biv}.
\epf

\section{Tightness -- Proof of Theorem~\ref{th_tight}}

In this section we will show that the sequence of random variables
$l_n(t) = n^{-1/2}L_n(t\sqrt n)$, $t\ge 0$, is tight
in C$[0,\infty)$. By \cite[p. 63]{KS} it suffices to prove tightness for C$[0,A]$. Hence we
consider $L_n(t)$ for $0\le t\le A\sqrt n$, where $A>0$ is an arbitrary real
constant.

By \cite[Theorem 12.3]{Bi68} tightness of $l_n(t) = n^{-1/2}L_n(t\sqrt n)$,
$0\le t\le A$, follows from tightness of $L_n(0)$ (which is trivial)
and from the existence of a constant $C>0$ such that \eqref{tight} holds 
for all non-negative integers $n,r,h$.

The fourth moment in Equation~\eqref{tight} can be expressed as the coefficient of
a suitable generating function. Indeed the generating function counting tree according to size as
well as the quantity $L_n(r)-L_n(r+h)$ is $\tilde y_{r,h}\(x,u,\rcp{u}\)$ ($\tilde y_{r,h}(x,u,v)$
is defined by \eqref{tightrec}) since assigning a weight $u$ to
the vertices in level $r$ and weight $1/u$ to those in level $r+h$ means that any tree having $n$
vertices and with $L_n(r)-L_n(r+h)=m$ contributes to the coefficient of $x^nu^m$ where $m$ can be
negative. The fourth moment can then be obtained by applying the operator 
$\(u\pdiff{}u{}\)^4$ and setting $u=1$ afterwards. Therefore we have 
$$
\E\left( L_n(r)-L_n(r+h)\right)^4 = \frac 1{y_n} [x^n] 
\U{\left( \pdiff{}u{}+7\pdiff{}u2+6\pdiff{}u3+\pdiff{}u4 \right)
\tilde y_{r,h}\(x,u,\rcp{u}\)}
$$
Thus, \eqref{tight} is equivalent to 
\begin{equation}
[x^n] 
\U{\left( \pdiff{}u{}+7\pdiff{}u2+6\pdiff{}u3+\pdiff{}u4 \right)
\tilde y_{r,h}\(x,u,\rcp{u}\)} \le C\,\frac{h^2}{\sqrt n} \rho^{-n}
\label{eqTh71}
\end{equation}
In order to prove \eqref{eqTh71} we use a result from 
\cite{FO90} saying that 
\[
F(x) = \Ord{(1-x/\rho)^{-\beta}}\qquad (x\in\Delta)
\]
implies
\[
[x^n] F(x) = \Ord{\rho^{-n} n^{\beta-1}},
\]
where $\Delta$ is the region of \eqref{eqDelta}

Hence, it is sufficient to show that
\begin{equation}\label{eqtoshow}
\U{\left( \pdiff{}u{}+7\pdiff{}u2+6\pdiff{}u3+\pdiff{}u4 \right)
\tilde y_{r,h}\(x,u,\rcp{u}\)} = \Ord{ \frac{h^2}{1-|y(x)|} }
= \Ord{ \frac{h^2}{\sqrt{|1-x/\rho|}} }
\end{equation}
for $x\in \Delta$ and $h\ge 1$. 
(Note that $\theta < \frac\pi 2$ implies that
$1 - |y(x)| \ge c \sqrt{|1-x/\rho|}$ for some constant $c>0$.)

Now we define
\[
\gamma_k^{[j]}(x) = \U{\frac{\partial^j y_k(x,u)}{\partial u^j} }
\qquad \mbox{and}\qquad
\gamma_{r,h}^{[j]}(x) = \U{\frac{\partial^j \tilde y_{r,h}\(x,u,\rcp{u}\)}
{\partial u^j} }
\]
and derive the following upper bounds.

\bl\label{Le8}
We have 
\begin{equation}\label{eqLe81}
\gamma_k^{[1]}(x) = \left\{ \begin{array}{ll}
\Ord{1} & \mbox{uniformly for $x\in \Delta$}, \\
 \Ord{|x/\rho|^k} & \mbox{uniformly for $|x|\le \rho$}
\end{array}\right.
\end{equation}
and 
\begin{equation}\label{eqLe82}
\gamma_{r,h}^{[1]}(x) = \left\{ \begin{array}{ll}
\Ord{\frac {h}{r+h}} & \mbox{uniformly for $x\in \Delta$}, \\
 \Ord{|x/\rho|^r} & \mbox{uniformly for $|x|\le \rho$,}
\end{array}\right.
\end{equation}
where $L$ is constant with $0<L< 1$.
\el

\bpf
We already know that $\gamma_k^{[1]}(x) = C_k(x) y(x)^k$,
where $C_k(x) = \Ord{1}$ and $|y(x)|\le 1$ for $x\in \Delta$.
Furthermore, by convexity we also have $|y(x)|\le |x/\rho|$ 
for $|x|\le \rho$. Hence, we obtain $\gamma_k^{[1]}(x) =
\Ord{|x/\rho|^k}$ for $|x|\le\rho$.

The functions $\gamma_{r,h}^{[1]}(x)$ are given by the recurrence
\[
\gamma_{r+1,h}^{[1]}(x) = y(x) \sum_{i\ge 1} \gamma_{r,h}^{[1]}(x^i)
\]
with initial value $\gamma_{0,h}^{[1]}(x) = y(x) - \gamma_h(x)$. To show this, just differentiate
\eqref{tightrec} w.r.t. $u$ and then plug in $u=1$. 
Hence, the representation $\gamma_{r,h}^{[1]}(x) = \gamma_r^{[1]}(x) -
 \gamma_{h+r}^{[1]}(x)$ follows by induction. Since, $\gamma_r^{[1]}(x)
 = (C(x)+O(L^r)) y(x)^r$ we thus get that
\[
\gamma_{r,h}^{[1]}(x) = \Ord{\sup_{x\in \Delta} |y(x)^r(1-y(x)^h)| + L^r}
\]
However, it is an easy exercise to show that 
\begin{equation}\label{eqOest}
\sup_{x\in \Delta} |y(x)^r(1-y(x)^h)| = \Ord{\frac h{r+h}}.
\end{equation}
For this purpose observe that if $x\in \Delta$ then we either
have $|y(x)-1|\le 1$ and $|y(x)|\le 1$, 
or $|y(x)|\le 1 - \eta$ for some $\eta> 0$.
In the second case we surely have 
\[
|y(x)^r(1-y(x)^h)| \le 2 (1-\eta)^r = \Ord{L^r}.
\]
For the first case we set $y = 1 - \rho e^{i\varphi}$ and observe
that 
\[
\left| 1- (1 - \rho e^{i\varphi})^h \right| \le (1+\rho)^h -1.
\]
Hence, if $r \ge 3h$ we thus obtain that
\[
|y(x)^r(1-y(x)^h)| \le \max_{0\le \rho\le 1}
(1-\rho)^r  \left( (1+\rho)^h -1 \right) \le \frac h{r-h} \le \frac {2h}{r+h}.
\]
If $r< 3h$ we obviously have
\[
|y(x)^r(1-y(x)^h)| \le 2 \le \frac {4h}{r+h}
\]
which completes the proof of (\ref{eqOest}).
Of course, we also have $L^r = \Ord{\frac h{h+r}}$. 
This completes the proof of the upper bound of $\gamma_{r,h}^{[1]}(x)$
for $x\in \Delta$.

Finally, the upper bound $\gamma_{r,h}^{[1]}(x) = \Ord{|x/\rho|^r}$
follows from (\ref{eqLe81}).
\epf

\bl\label{Le9}
We have 
\begin{equation}\label{eqLe91}
\gamma_k^{[2]}(x) = \left\{ \begin{array}{ll}
\Ord{\min\left\{k, \frac 1{1-|y(x)|}\right\}} & 
\mbox{uniformly for $x\in \Delta$}, \\
 \Ord{|x/\rho|^k} & \mbox{uniformly for $|x|\le \rho-\eta$}
\end{array}\right.
\end{equation}
and 
\begin{equation}\label{eqLe92}
\gamma_{r,h}^{[2]}(x) = \left\{ \begin{array}{ll}
\Ord{\min\left\{h, \frac 1{1-|y(x)|}\right\}} & 
\mbox{uniformly for $x\in \Delta$}, \\
 \Ord{|x/\rho|^r} & \mbox{uniformly for $|x|\le \rho-\eta$}
\end{array}\right.
\end{equation}
for every $\eta>0$.
\el

\brem
By doing a more precise analysis similarly to Lemma~\ref{Le7}
we can, for example, show that $\gamma_k^{[2]}(x)$ can be represented as
\begin{equation}\label{eqgamma2rep}
\gamma_k^{[2]}(x) = y(x)^k \sum_{\ell = 1}^k D_{k,\ell}(x) y(x)^{\ell -1},
\end{equation}
where the functions $D_{k,\ell}(x)$ are analytic in $\Delta$.
For every $\ell$ there is a limit $D_{\ell}(x) = \lim_{k\to\infty} D_{k,\ell}(x)$
with
\[
D_{k,\ell}(x) = D_\ell(x) + O(\tilde L^{k+\ell}),
\]
where $0< \tilde L < 1$.
Furthermore these limit functions $D_\ell(x)$ satisfy
\[
D_\ell(x) = C(x)^2 + O(\tilde L^\ell).
\]
Since we will not make use of this precise representation
we leave the details to the reader.
\erem

\bpf
The bound $\gamma_k^{[2]}(x) = O(|x/\rho|^k)$ (for $|x|\le \rho-\eta$)
and the bound $\gamma_k^{[2]}(x) = O(k)$ follow from 
Lemma~\ref{Le7-2}. In order to complete the analysis for
 $\gamma_{k}^{[2]}(x)$ we recall the recurrence derived from \eqref{rec} (compare also with
\eqref{eqgammak2rec-0})
\begin{equation}\label{eqgammak2rec}
\gamma_{k+1}^{[2]}(x) = y(x) \sum_{i\ge 1} i\gamma_k^{[2]}(x^i) 
+ y(x) \left(\sum_{i\ge 1} \gamma_k^{[1]}(x^i) \right)^2  + 
y(x) \sum_{i\ge 2} (i-1)\gamma_k^{[1]}(x^i) 
\end{equation}
that we rewrite to 
\begin{equation} \label{bkrec}
\gamma_{k+1}^{[2]}(x) = y(x) \gamma_{k}^{[2]}(x) + b_k(x),
\end{equation} 
where
\begin{equation} \label{bkdef}
b_k(x) = y(x) \sum_{i\ge 2} i\gamma_k^{[2]}(x^i) 
+ y(x) \left(\sum_{i\ge 1} \gamma_k^{[1]}(x^i) \right)^2  + 
y(x) \sum_{i\ge 2} (i-1)\gamma_k^{[1]}(x^i).
\end{equation} 
Note that for $i\ge 2$ we have $|x^i|\le\rho-\eta$. Therefore we can apply the second estimate of
\eqref{eqLe91} (first and third sum of \eqref{bkdef}) and the estimate \eqref{eqLe81} (the first
for $i=1$ in the second sum of \eqref{bkdef}, the second for the other summands) and obtain then 
\begin{equation} \label{bkbound}
b_k(x) = \Ord{1} \mbox{ uniformly for
$x\in \Delta$.}
\end{equation}  

Since $\gamma_{0}^{[2]}(x) = 0$, the solution of the recurrence \eqref{bkrec} 
can be written as
\begin{equation} \label{solution}
\gamma_{k}^{[2]}(x) = b_{k-1}(x) + y(x) b_{k-2}(x) + \cdots + y(x)^{k-1} b_0(x).
\end{equation} 
So \eqref{bkbound} implies finally 
\[
\gamma_{k}^{[2]}(x) = \Ord{\sum_{j=0}^{k-1} |y(x)|^j} = \Ord{\frac 1{1-|y(x)|} }.
\]
which completes the proof of (\ref{eqLe91}).

The recurrence for $\gamma_{r,h}^{[2]}(x)$ is similar to that
of $\gamma_{k}^{[2]}(x)$:
\begin{equation}\label{eqgammarh2rec}
\gamma_{r+1,h}^{[2]}(x) = y(x) \sum_{i\ge 1} i\gamma_{r,h}^{[2]}(x^i) 
+ y(x) \left( \sum_{i\ge 1} \gamma_{r,h}^{[1]}(x^i)\right)^2 + 
y(x) \sum_{i\ge 2} (i-1)\gamma_{r,h}^{[1]}(x^i) 
\end{equation}
with initial value $\gamma_{0,h}^{[2]}(x) = \gamma_{h}^{[2]}(x)$.
We again use induction.
Assume that we already know that $|\gamma_{r,h}^{[2]}(x)| \le D_{r,h} |x/\rho|^k$
for $|x|\le \rho-\eta$ and for some constant $D_{r,h}$.
By (\ref{eqLe92}) we can set $D_{0,h} = D_h$ which
is bounded as $h\to\infty$.
We also assume that 
$|\gamma_{r,h}^{[1]}(x)| \le C |x/\rho|^k$ for $|x|\le \rho-\eta$.
Then by (\ref{eqgammarh2rec}) we get
\begin{align*}
|\gamma_{r+1,h}^{[2]}(x)| 
&\le  D_{r,h} |x/\rho|^{k+1} + D_{r,h} |x/\rho|
\frac{2 |x|^{2k}/\rho^k}{(1-|x|^k)^2} \\
&+ C^2 |x/\rho| \left( \frac{|x/\rho|^{k}}{1-|x|^k}\right)^2 + 
C |x/\rho| \frac{2 |x|^{2k}/\rho^k}{(1-|x|^k)^2} .
\end{align*}
Thus, we can set
\[
D_{r+1,h} = D_{r,h} \left( 1+ \frac{2 (\rho-\eta)^{k}}{(1-\rho^k)^2} \right) + 
C^2 \frac{(\rho-\eta)^{k}}{(1-\rho^k)^2} + 
C\frac{2 (\rho-\eta)^{k}}{(1-\rho^k)^2} 
\]
which shows that the constants $D_{r,h}$ are uniformly bounded. Consequently
$\gamma_{r,h}^{[2]}(x) = \Ord{|x/\rho|^k}$ for $|x|\le \rho-\eta$.

Next we start from \eqref{eqgammarh2rec} and 
 assume that $|\gamma_{r,h}^{[2]}(x)| \le \bar D_{r,h}$
for $x\in \Delta$. We already know that
$|\gamma_{r,h}^{[1]}(x)| \le C \frac h{h+k}$ for $x\in \Delta$.
Hence,
\begin{align*}
|\gamma_{r+1,h}^{[2]}(x)| &\le \bar D_{r,h} + D_{r,h} 
\sum_{i\ge 2} i |x^i/\rho|^k \\
&+ C^2 \left( \frac h{k+h} + \sum_{i\ge 2} |x^i/\rho|^k\right)^2 + 
C \sum_{i\ge 2} (i-1)|x^i/\rho|^k\\
&\le \bar D_{r,h} + 8 D_{r,h}  (\rho+\eta)^{2k} /\rho^k \\
&+ C^2\left( \frac h{k+h} + 2(\rho+\eta)^{2k}/\rho^k\right)^2 
+ 4C  (\rho+\eta)^{2k} /\rho^k.
\end{align*}
Thus, we can set
\begin{align*}
\bar D_{r+1,h} &=  \bar D_{r,h} 
+ 8 D_{r,h}  (\rho+\eta)^{2k} /\rho^k 
+ C^2\left( \frac h{k+h} + 2(\rho+\eta)^{2k}/\rho^k\right)^2 \\
&+ 4C  (\rho+\eta)^{2k} /\rho^k
\end{align*}
with initial value $\bar D_{0,h} = \bar D_h = \Ord{h}$ and
obtain a uniform upper bound of the form
\[
\bar D_{r,h} = \Ord{h}.
\]
Consequently $\gamma_{r,h}^{[2]}(x)= \Ord{h}$ for $x\in \Delta$.

Thus, in order to complete the proof of (\ref{eqLe92}) it remains to 
prove $\gamma_{r,h}^{[2]}(x) = \Ord{1/(1-|y(x)|)}$ 
for $x\in \Delta$. Analogously to the way we obtained \eqref{solution} from \eqref{eqgammak2rec}
we obtain from \eqref{eqgammarh2rec} the representation of $\gamma_{r,h}^{[2]}(x)$ as
\begin{equation}\label{eqgrh2rep}
\gamma_{r,h}^{[2]}(x) = \gamma_{0,h}^{[2]}(x) + 
c_{k-1,h}(x) + y(x) c_{k-2,h}(x) + \cdots + y(x)^{k-1} c_{0,h}(x),
\end{equation}
where 
\[
c_{j,h}(x) = y(x) \sum_{i\ge 2} i\gamma_{j,h}^{[2]}(x^i) 
+ y(x) \left( \sum_{i\ge 1} \gamma_{j,h}^{[1]}(x^i)\right)^2 + 
y(x) \sum_{i\ge 2} (i-1)\gamma_{j,h}^{[1]}(x^i).
\]
Observe that there exists $\eta>0$ such that
$|x^i|\le \rho-\eta$ for $i\ge 2$ and $x\in \Delta$. Hence it
follows in a similar fashion as we showed \eqref{bkbound} 
that $c_{j,h}(x) = \Ord{1}$ for $x\in \Delta$.
Since $\gamma_{0,h}^{[2]}(x) = \gamma_{h}^{[2]}(x) = \Ord{1/(1-|y(x)|)}$,
we consequently get
\[
\gamma_{r,h}^{[2]}(x) = \gamma_{h}^{[2]}(x)+ \Ord{\frac 1{1-|y(x)|}}
= \Ord{\frac 1{1-|y(x)|}}.
\]
\epf

\brem
Note that the theorem our proof relies on, namely \cite[Theorem 12.3]{Bi68}, actually requires the
existence of $\alpha>0$ and $\beta>1$ such that 
$$
{\bf E}\, |L_n(r) - L_n(r+h)|^\alpha = \Ord{ (h \sqrt n)^\beta}. 
$$ 
The estimates of Lemma~\ref{Le9} already prove that
\[
{\bf E}\, (L_n(r) - L_n(r+h))^2 = \Ord{h \sqrt n}.
\]
Unfortunately this estimate is slightly too weak to prove tightness. Third moments are technically
unpleasant they attain positive and negative signs. So we actually 
have to deal with 4-th moments. 
\erem

Before we start with bounds for $\gamma_k^{[3]}(x)$ and 
$\gamma_k^{[4]}(x)$ we need an auxiliary bound.

\bl\label{Le9-2}
We have uniformly for $x\in \Delta$
\begin{equation}\label{eqLe9-21}
\sum_{r\ge 0} |\gamma_{r,h}^{[1]}(x)\gamma_{r,h}^{[2]}(x)| = 
\Ord{h^2}.
\end{equation}
\el

\bpf
We use the representation (\ref{eqgrh2rep}), where we
can approximate $c_{j,h}(x)$ by 
\[
c_{j,h}(x) = y(x) \gamma_{j,h}^{[1]}(x)^2 + \Ord{L^j}= 
\Ord{ \frac{h^2}{(r+h)^2} }
\] 
uniformly for $x\in \Delta$ with some constant $L$ that satisfies
$0< L< 1$. Furthermore, we use the approximation
\[
\gamma_{r,h}^{[1]}(x) = C(x)y(x)^r(1-y(x)^h) + \Ord{L^r} 
\]
that is uniform for $x\in \Delta$. For example, this shows
\[
\sum_{r\ge 0} |\gamma_{r,h}^{[1]}(x)| = 
|C(x)|\frac{|1-y(x)^h|}{1-|y(x)|} + \Ord{1}.
\]
Now observe that for $x\in \Delta$ there exists a constant $c> 0$
with $|1-y(x)| \ge c (1 - |y(x)|)$. Hence it follows that
\[
\frac{|1-y(x)^h|}{1-|y(x)|} = \Ord{\left| \frac{1-y(x)^h}{1-y(x)}\right|} = 
\Ord{h}
\]
and consequently 
\[
\sum_{r\ge 0} |\gamma_{r,h}^{[1]}(x)| = \Ord{h}.
\]
Similarly we get
\begin{align*}
\sum_{r\ge 1} &|\gamma_{r,h}^{[1]}(x)|
\, \left| \sum_{j< r} y(x)^{r-j-1}  c_{j,h}(x) \right|
\le \sum_{j\ge 0} |c_{j,h}(x)|\, |y(x)|^{-j-1} 
\sum_{r>j} |y(x)|^r  |\gamma_{r,h}^{[1]}(x)|\\
&= \sum_{j\ge 0} |c_{j,h}(x)|\, |y(x)|^{-j-1} 
\left( |C(x)| \, |y(x)|^{2j+2} \frac{|1-y(x)^h|}{1-|y(x)|^2} 
+ \Ord{|y(x)|^{j+1}L^j}\right) \\
&= \Ord{\sum_{j\ge 0} \frac{h^3}{(j+h)^2}} \\
&= \Ord{h^2}.
\end{align*}
Hence, we finally obtain  
\begin{align*}
\sum_{r\ge 0} |\gamma_{r,h}^{[1]}(x)\gamma_{r,h}^{[2]}(x)| &\le
 \sum_{r\ge 0} |\gamma_{r,h}^{[1]}(x)| 
\, |\gamma_{0,h}^{[2]}(x)|
+  \sum_{r\ge 1} |\gamma_{r,h}^{[1]}(x)|\,
\left| \sum_{j< r} y(x)^{r-j-1}  c_{j,h}(x)
\right| \\
&=|\gamma_{h}^{[2]}(x)|\cdot\sum_{r\ge 0} |\gamma_{r,h}^{[1]}(x)|
+  \sum_{r\ge 1} |\gamma_{r,h}^{[1]}(x)|\,
\left| \sum_{j< r} y(x)^{r-j-1}  c_{j,h}(x)
\right|
&= \Ord{ h^2}.
\end{align*}
where we used the fact $|\gamma_{h}^{[2]}(x)|=\Ord{h}$ from Lemma~\ref{Le9}.
\epf

\bl\label{Le10}
We have 
\begin{equation}\label{eqLe101}
\gamma_k^{[3]}(x) = \left\{ \begin{array}{ll}
\Ord{\min \{ k^2, \frac k{1-|y(x)|}\} } & \mbox{uniformly for $x\in \Delta$}, \\
 \Ord{|x/\rho|^k} & \mbox{uniformly for $|x|\le \rho-\eta$}
\end{array}\right.
\end{equation}
and 
\begin{equation}\label{eqLe102}
\gamma_{r,h}^{[3]}(x) = \left\{ \begin{array}{ll}
\Ord{ \min \{ h^2,\frac h{1-|y(x)|} \} }  & \mbox{uniformly for $x\in \Delta$}, \\
 \Ord{|x/\rho|^r} & \mbox{uniformly for $|x|\le \rho-\eta$}
\end{array}\right.
\end{equation}
for every $\eta>0$.
\el

\bpf
The recurrence for $\gamma_k^{[3]}(x)$ is again obtain by differentiation of \eqref{rec} and 
given by
\begin{align}
\gamma_{k+1}^{[3]}(x) &= y(x) \sum_{i\ge 1} i^3\gamma_k^{[3]}(x^i) +
y(x) \left( \sum_{i\ge i} \gamma_k^{[1]}(x^i) \right)^3 
+ 3 y(x) \left( \sum_{i\ge 1}  \gamma_k^{[1]}(x^i) \right)
\left( \sum_{i\ge 1} i \gamma_k^{[2]}(x^i) \right) \nonumber \\
&+ 3y(x)  \left( \sum_{i\ge 1}  \gamma_k^{[1]}(x^i) \right)
\left( \sum_{i\ge 1} (i-1) \gamma_k^{[i]}(x^i) \right) 
+ 3y(x) \sum_{i\ge 1} i(i-1)\gamma_k^{[2]}(x^i) \\
&+ y(x) \sum_{i\ge 1} (i-1)(i-2)\gamma_k^{[1]}(x^i) \nonumber
\end{align}
By inspecting the proof of Lemmas~\ref{Le8} and \ref{Le9} one expects
that the only {\it important} part of this recurrence if given by
\begin{equation}\label{eqrecgk3simple}
\gamma_{k+1}^{[3]}(x) = y(x)\gamma_k^{[3]}(x) +
y(x) \gamma_k^{[1]}(x)^3 
+ 3 y(x) \gamma_k^{[1]}(x)  \gamma_k^{[2]}(x) + R_k
\end{equation}
and $R_k$ collects the {\it less important} remainder terms that
only contributes exponentially small terms. Thus, in order to
shorten our presentation we will only focus on these terms.
In particular it is easy to show the bound 
$\gamma_{k}^{[3]}(x) = \Ord{|x/\rho|^k}$
for $|x|\le \rho-\eta$. (We omit the details.)

Next, since $y(x) \gamma_k^{[1]}(x)^3 
+ 3 y(x) \gamma_k^{[1]}(x)  \gamma_k^{[2]}(x) + R_k = \Ord{k}$ for $x\in\Delta$,   
it directly follows that $\gamma_k^{[3]}(x) = \Ord{k^2}$.

Now we proceed by induction and observe that a bound of the 
form $|\gamma_{k}^{[3]}(x)| \le E_k /(1-|y(x)|)$ leads to
\begin{align*}
|\gamma_{k+1}^{[3]}(x)| &\le \frac {E_k}{1-|y(x)|} + 
\Ord{\frac 1{1-|y(x)|}} + |R_k|
\end{align*}
and consequently  to $E_{k+1} \le E_k + \Ord{1}$. 
Hence, $E_k = \Ord{k}$ and  $\gamma_k^{[3]}(x) = 
\Ord{k/(1-|y(x)|)}$.

Similarly, the leading part of the recurrence for 
$\gamma_{r,h}^{[3]}(x)$ is given by
\begin{align}\label{eqrecgrh3simple}
\gamma_{r+1,h}^{[3]}(x) &= y(x)\gamma_{r,h}^{[3]}(x) +
y(x) \gamma_{r,h}^{[1]}(x)^3 
+ 3 y(x) \gamma_{r,h}^{[1]}(x)  \gamma_{r,h}^{[2]}(x) + \bar R_{r,h} \\
&= y(x)\gamma_{r,h}^{[3]}(x) + d_{r,h}(x), \nonumber
\end{align}
where 
\[
d_{r,h}(x) = y(x) \gamma_{r,h}^{[1]}(x)^3 + 3 y(x) \gamma_{r,h}^{[1]}(x)  \gamma_{r,h}^{[2]}(x) + \bar R_{r,h} = \Ord{h}
\]
and the initial value is given by
\[
\gamma_{0,h}^{[3]}(x) = - \gamma_{h}^{[3]}(x)
- 3 \gamma_{h}^{[2]}(x) = \Ord{\min\left\{ h^2, \frac h{1-|y(x)|}   \right\}}.
\]
Note that we also assume that 
$\gamma_{r,h}^{[3]}(x) = \Ord{|x/\rho|^r}$
for $|x|\le \rho-\eta$ (which can be easily proved).
Consequently it directly follows that
\begin{align*}
\gamma_{r,h}^{[3]}(x) &= \gamma_{0,h}^{[3]}(x) 
+ d_{r-1,h}(x) + y(x) d_{r-1,h}(x) + \cdots + y(x)^{r-1} d_{0,h}(x)\\
&= \Ord{ \frac h{1-|y(x)|} }.
\end{align*}
Next observe that Lemmas~\ref{Le8}--\ref{Le9-2} ensure that
\[
\sum_{j\ge 0} |d_{j,h}(x)| = \Ord{h^2}
\]
uniformly for $x\in \Delta$. Hence, we finally get
\[
\gamma_{r,h}^{[3]}(x) = \Ord{h^2}
\]
which completes the proof of Lemma~\ref{Le10}.
\epf

\bl\label{Le11}
We have 
\begin{equation}\label{eqLe111}
\gamma_k^{[4]}(x) = \left\{ \begin{array}{ll}
\Ord{\frac {k^2}{1-|y(x)|}} & \mbox{uniformly for $x\in \Delta$}, \\
 \Ord{|x/\rho|^k} & \mbox{uniformly for $|x|\le \rho-\eta$}
\end{array}\right.
\end{equation}
and 
\begin{equation}\label{eqLe112}
\gamma_{r,h}^{[4]}(x) = \left\{ \begin{array}{ll}
\Ord{\frac {h^2}{1-|y(x)|}} & \mbox{uniformly for $x\in \Delta$}, \\
 \Ord{|x/\rho|^r} & \mbox{uniformly for $|x|\le \rho-\eta$}
\end{array}\right.
\end{equation}
for every $\eta>0$.
\el

\bpf
The proof is very similar to that of Lemma~\ref{Le10}.
First, the recurrence for $\gamma_k^{[4]}(x)$ is essentially
of the form
\begin{align}\label{eqrecgk4simple}
\gamma_{k+1}^{[4]}(x) &= y(x)\gamma_k^{[4]}(x) +
y(x) \gamma_k^{[1]}(x)^4 
+ 4 y(x) \gamma_k^{[1]}(x)  \gamma_k^{[3]}(x) \\
&+ 6 y(x) \gamma_k^{[1]}(x)^2  \gamma_k^{[2]}(x)
+ 3 y(x)   \gamma_k^{[2]}(x)^2 +
R_k,   \nonumber
\end{align}
where $R_k$ collects all exponentially small summands.
We assume that we have already proved 
the upper bound $\gamma_{k}^{[4]}(x) = \Ord{|x/\rho|^k}$
for $|x|\le \rho-\eta$. Now, by induction and 
the assumption $|\gamma_{k}^{[4]}(x)| \le F_k /(1-|y(x)|)$ and
the known estimates $\gamma_{k}^{[1]}(x) = \Ord{1}$, 
$\gamma_{k}^{[2]}(x) = \Ord{\min\{k, 1/(1-|y(x)|)\}}$, and
$\gamma_{k}^{[3]}(x) = \Ord{k/(1-|y(x)|)}$ we get
\begin{align*}
|\gamma_{k+1}^{[4]}(x)| &\le \frac {F_k}{1-|y(x)|} + \Ord{|y(x)|^k} + 
\Ord{\frac k{1-|y(x)|}}+ \Ord{\frac 1{1-|y(x)|}}+ |R_k|
\end{align*}
and consequently $F_k = \Ord{k^2}$.

Finally, the essential part of the recurrence for 
$\gamma_{r,h}^{[4]}(x)$ is given by
\begin{align}\label{eqrecgrh4simple}
\gamma_{r+1,h}^{[4]}(x) &= y(x)\gamma_{r,h}^{[4]}(x) +
y(x) \gamma_{r,h}^{[1]}(x)^4 
+ 4 y(x) \gamma_{r,h}^{[1]}(x)  \gamma_{r,h}^{[3]}(x) \\
&+ 6 y(x) \gamma_{r,h}^{[1]}(x)^2  \gamma_{r,h}^{[2]}(x)
+ 3 y(x)   \gamma_{r,h}^{[2]}(x)^2 + \bar R_{r,h}    \nonumber \\
&= y(x) \gamma_{r,h}^{[4]}(x)  + e_{r,h}(x), \nonumber
\end{align}
where 
\begin{align*}
e_{r,h}(x) &= y(x) \gamma_{r,h}^{[1]}(x)^4 
+ 4 y(x) \gamma_{r,h}^{[1]}(x)  \gamma_{r,h}^{[3]}(x) \\
&+ 6 y(x) \gamma_{r,h}^{[1]}(x)^2  \gamma_{r,h}^{[2]}(x)
+ 3 y(x)   \gamma_{r,h}^{[2]}(x)^2 + \bar R_{r,h}.
\end{align*}
As above, $\bar R_{r,h}$ collects all exponentially small terms.
Thus,
\[
\gamma_{r,h}^{[4]}(x) =   \gamma_{0,h}^{[4]}(x) 
+ e_{r-1,h}(x) + y(x) e_{r-1,h}(x) + \cdots + y(x)^{r-1} e_{0,h}(x).
\]
If we use the known estimates $\gamma_{r,h}^{[1]}(x) = \Ord{1}$, 
$\gamma_{r,h}^{[2]}(x) = \Ord{h }$, and
$\gamma_{r,h}^{[3]}(x) = \Ord{h^2 }$ which gives $d_{r,h} = \Ord{h^2}$ 
and the initial condition
\[
\gamma_{0,h}^{[4]}(x) = 12 \gamma_{h}^{[2]}(x) + 
8 \gamma_{h}^{[3]}(x) + \gamma_{h}^{[4]}(x) = \Ord{\frac {h^2}{1-|y(x)|}}
\]
we obtain
\[
\gamma_{r,h}^{[4]}(x)  = \Ord{\frac {h^2}{1-|y(x)|}}.
\]
This completes the proof of Lemma~\ref{Le11}.
\epf

The proof of (\ref{eqtoshow}) is now immediate. As already 
noted this implies (\ref{eqTh71}) and proves Theorem~\ref{th_tight}.

\section{The Height -- Proof of Theorem~\ref{height}}

Let $y_n^{(k)}$ denote the number of trees with $n$ nodes and height less than $k$. Then 
the generating function $y_k(x)=\sum_{n\ge 1} y_n^{(k)}x^n$ satisfies the recurrence relation
\begin{align*}
y_0(x)&=0 \nonumber \\
y_{k+1}(x)&=x\exp\(\sum_{i\ge 1}\frac{y_k(x^i)}i\),\quad k\ge 0.
\end{align*}
Obviously $y_k(x)=y_k(x,0)$ where the function on the right-hand side is the generating function
of \eqref{rec} which we used to analyze the profile in the previous sections. So $w_k$ and
$\Sigma_k$ could be defined accordingly. However, the proof given here relies heavily on the
seminal work of Flajolet and Odlyzko \cite{FO82} on the height of binary trees. Therefore, to be
in accordance with the notation used there, we work with the opposite sign and set 
\[
e_k(x) = y(x)- y_k(x), 
\]
that is $e_k(x) = - w_k(x,0)$. Then $e_k$ satisfies the recurrence
\begin{equation}\label{eqrecek}
e_{k+1}(x) = y(x)\left( 1 - e^{-e_k(x)-E_k(x)} \right),
\end{equation}
where
\begin{equation} \label{Ekdef}
E_k(x) = \sum_{i\ge 2} \frac {e_k(x^i)}{i} = - \Sigma_k(x,0).
\end{equation} 
The function $e_k(x)$ is the generating function for the number of trees with height at least $k$.

The proof of Theorem~\ref{height} follows the same principles as the proof of the
corresponding properties of the height of Galton-Watson trees
(see \cite{FO82,FGOR93}). However, the term $E_k(x)$ 
needs some additional considerations. 

\bpr 
If $e_k(x)$ satisfies 
\begin{equation}\label{eqekrep}
e_k(x) = \frac{y(x)^{k}}
{\frac 12 \frac{1-y(x)^k}{1-y(x)} +
O\left( \min \left\{ \log k, \log \frac 1{1-|y(x)|} \right\} \right)}.
\end{equation}
for $x\in\Delta_\eps$, then Theorem~\ref{height} follows.  
\epr 

\bpf
Since $e_k(x)$ is bounded in $\Delta\setminus\Delta_\eps$, only the local behaviour near $x=\rho$
determines the asymptotic height of P\'olya trees. 
The shape \eqref{eqekrep} of $e_k(x)$ precisely matches that of the
corresponding quantity for simply generated trees. Flajolet and Odlyzko showed that
\eqref{eqekrep} implies \eqref{av_height}, see \cite[p.~204]{FO82} where this argument was used to 
derive the average height as well as the other moments of the height of simply generated trees. 
\epf

\bpr\label{keyprop}
Suppose that for $x\in\Delta_\eps$ the estimate 
\begin{equation}\label{eqProOLk}
\frac{|E_k(x)|}{|e_k(x)|^2} = O(L^k)
\end{equation}
holds for some $L< 1$ and that 
\begin{equation} \label{einsdurchk} 
|e_k(x)y(x)^k|=\Ord{1/k}.
\end{equation}  
Then we have \eqref{eqekrep} for $x\in\Delta_\eps$.
\epr

\bpf
Equation \eqref{eqrecek} can be rewritten to (omitting the argument $x$)
\[
e_{k+1} = y e_k \left( 1 -\frac {e_k}2 + 
O\left( e_k^2 + \frac{E_k}{e_k} \right) \right),
\]
resp.\ to
\[
\frac y{e_{k+1}} = \frac 1{e_k} + \frac 12 + 
O\left( e_k +  \frac{E_k}{e_k^2} \right).
\]
This leads to the representation
\begin{equation}\label{eqtkek}
\frac{y^k}{e_k} = \frac 1{e_0} + \frac 12 \frac{1-y^k}{1-y} 
+ O\left( \sum_{\ell < k} |e_\ell y^\ell| \right) + 
O\left( \sum_{\ell < k} \frac{|E_\ell|}{|e_\ell^2|} |y^\ell| \right).
\end{equation}
Recall that $e_0=y(x)$. By (\ref{eqProOLk}) and $\rcp{|y(x)|}=\Ord{1}$ (for $x\in\Delta_\eps$), 
a consequence of Lemma~\ref{lem_2}, this implies 
\begin{equation} \label{zugrob}
e_k = \frac{y^k}{\frac 12 \frac{1-y^k}{1-y} 
+ O\left( \sum_{\ell < k} |e_\ell y^\ell| \right) + O(1)}
\end{equation} 
and \eqref{einsdurchk} yields \eqref{eqekrep} after all.
\epf 

\brem
Note that the proof above does not make explicit use of the domain of $x$. Thus the implication of
Proposition~\ref{keyprop} is still true if we write, for instance,  $0\le x\le \rho$ instead of
$x\in\Delta_\eps$ in \eqref{eqekrep}, \eqref{eqProOLk}, and \eqref{einsdurchk}. We remark that we 
will use such modifications of Proposition~\ref{keyprop} in the sequel, though we do not state 
several almost identical propositions differing only in the domain of $x$. 
\erem

Note that (\ref{eqtkek}) and (\ref{eqekrep}) can be made more precise.
Set 
\[
S_k = \frac{e_k^2 \left(e^{-E_k} -1 \right)}
{\left( e^{e_k}-1 \right) \left( 1 -e^{-e_k-E_k}\right) },
\]
and define a function $h(v)$ by
\[
\frac v{1 - e^{-v}} = 1 + \frac v2 + v^2 h(v).
\]
Then the recurrence $e_{k+1} = y(1- e^{-e_k-E_k})$ rewrites to
\[
\frac y{e_{k+1}} = \frac 1{e_k} + \frac 12 + e_k h(e_k) + 
\frac {S_k}{e_k^2}
\]
and leads to the explicit representations
\begin{equation}\label{eqekexprep1}
\frac{y^k}{e_k} = \frac 1{e_0} + \frac 12 \frac{1-y^k}{1-y}
+ \sum_{\ell < k} e_\ell h(e_\ell) y^\ell + 
\sum_{\ell < k} \frac{S_\ell}{e_\ell^2} y^\ell
\end{equation}
and
\begin{equation}\label{eqekexprep2}
e_k = \frac{y^k}{\frac 1{e_0} + \frac 12 \frac{1-y^k}{1-y}
+ \sum_{\ell < k} e_\ell h(e_\ell) y^\ell + 
\sum_{\ell < k} \frac{S_\ell}{e_\ell^2} y^\ell  }.
\end{equation}
This formula is a refinement of \eqref{eqekrep} since it makes the error term explicit. We will
use it in the sequel. 

Furthermore, note that if we just assume $e_k\to 0$ and $E_k = o(e_k)$ as 
$k\to\infty$, then 
\begin{equation} \label{auxexpr}
S_k \sim - E_k \mbox{ and } h(e_k)=\Ord{1}.
\end{equation}

We start our precise analysis with an a priori bound for $e_k(x)$. The next step is proving
\eqref{eqekrep} for $0\le x\le \rho$. Then we will, little by little, enlarge the allowed domain
for $x$ and arrive finally at Proposition~\ref{keyprop} as stated above.

\bl\label{upperbound_height}
Let $|x|\le\rho$. Then there is a $C>0$ such that 
$$
|e_k(x)|\le\frac{C}{\sqrt k}\btr{\frac x\rho}^k.
$$
\el

\bpf
Obviously, we have 
$$
|e_k(x)|=\sum_{n>k} (y_n-y_n^{(k)}) |x|^n \le \sum_{n>k} y_n |x|^n. 
$$
The assertion follows now from  $y_n\sim c \rho^{-n} n^{-3/2}$
for some constant $c>0$.
\epf

Lemma~\ref{upperbound_height} applies to $E_k(x)$.
\bcor\label{CorEk}
Suppose that $|x|< \sqrt \rho$. Then there exists a constant $C_0>0$
with
\[
|E_k(x)| \le \frac{C_0}{\sqrt k} \left| \frac {x^2}{\rho} \right|^k.
\]
\ecor

\brem
Observe that in the definition of $E_k(x)$, Equation~\eqref{Ekdef},  
the arguments in the sum are raised to a power of at
least 2. Therefore $E_k(x)$ is an analytic function in the domain $|x|< \sqrt \rho$ and not only
in the smaller domain $|x|\le \rho$.
\erem

The next lemma shows that $e_k(x)$ behaves as expected if $x$ is
on the positive real axis.
\bl\label{Lerecreal}
Suppose that $0\le x \le \rho$ is real. Then (\ref{eqekrep}) holds in this domain for $x$. 
\el

\brem
As remarked before, we will show a weaker version of \eqref{eqekrep}. During the proof we will
make use of the analogous weaker versions of \eqref{eqProOLk} and \eqref{einsdurchk}. Therefore,
in the following proof a reference to one of these formulas means the referred formula, but 
for $0\le x\le\rho$ and not for $x\in\Delta_\eps$. 
\erem

\bpf
Let $\tilde e_k(x)$ be defined by $\tilde e_0(x) = y(x)$ and by
$\tilde e_{k+1}(x) = y(x)(1 - e^{-\tilde e_k(x)})$ (for $k \ge 0$).
Then $\tilde e_k(x)$ is precisely the analogue of $e_k$ for Cayley trees, a class of simply
generated trees (preceisely: the class of labelled rooted trees). 
So $\tilde e_k(x)$ behaves like \eqref{eqekrep} in $\Delta$. 

However, if $0\le x\le \rho$ then we obtain by induction that 
$e_k(x)\ge \tilde e_k(x)$. Hence, by combining (\ref{eqekrep})
with the upper bound from Lemma~\ref{upperbound_height} we have
\[
\frac{E_k(x)}{e_k(x)^2} \le \frac{E_k(x)}{\tilde e_k(x)^2} = O(L^k)
\]
for some $L$ with $0<L< 1$. 
Thus (\ref{eqProOLk}) is satisfied.

In order to show the second assumption \eqref{einsdurchk} 
of Proposition~\ref{keyprop} note that by Lemma~\ref{upperbound_height} $e_k(x)$ is even
exponentially small for $x<\rho$. For the case $x=\rho$ observe that \eqref{eqProOLk} in
conjunction with Lemma~\ref{upperbound_height} guarantee \eqref{auxexpr}. Applying this to 
\eqref{eqekexprep2} implies 
$$
e_k(x) = \frac{y(x)^{k}}
{\frac 12 \frac{1-y(x)^k}{1-y(x)} +
\Ord{\sqrt k}}.
$$
This equation yields 
\begin{equation} \label{ek2k} 
e_k(\rho)\sim \frac 2k
\end{equation} 
and completes the proof. 
\epf

The analysis of $e_k(x)$ for complex $x$ with $|x|\le \rho$ is 
not too difficult. The next two lemmas consider the case
$|x|\le \rho$ and $|x-\rho|\le \varepsilon$ and the case
$|x|\le \rho-\varepsilon$.

\bl\label{Leek1}
There exists $\varepsilon>0$ such that (\ref{eqekrep})
holds for all $x$ with $|x|\le \rho$ and $|x-\rho|\le \varepsilon$.
\el
\bpf
First recall that $|e_k(x)| \le C/\sqrt k$ and $|y(x)|<1$ for $x\neq\rho$. Moreover, in the proof of
the previous lemma we showed $e_k(\rho) \sim 2/k$. Hence \eqref{einsdurchk} is satisfied. 

Suppose that we can show that $|E_k/e_k^2|\le 1$ or, equivalently, $|E_k/e_k|\le |e_k|\le C/\sqrt
k$. Then it follows that 
\begin{align}
|e_{k+1}|&=|y|\, |e_k|\, \left|\frac{1-e^{-e_k-E_k}}{e_k} \right| \nonumber \\
&=|y|\, |e_k|\, \(1+\frac{E_k}{e_k}+\Ord{\frac{(e_k+E_k)^2}{e_k}}\) \nonumber \\
&\ge |y|\, |e_k|\, 1-C_1|e_k| \nonumber \\
&\ge  |y|\, |e_k| e^{-C_2 k^{-1/2}}. \label{eqek1ek}
\end{align}
where $C_1,C_2$ are suitable constants.

Now we choose $k_0$ sufficiently large such that 
\[
e^{-2C_2 \sqrt{k}} \le \frac 1{k} \quad \mbox{and}\quad
C_0 \rho^{k/2} e^{4C_2 \sqrt{k}} \le 1
\]
hold for all $k\ge k_0$.
By continuity, \eqref{ek2k} implies the existence of an $\varepsilon> 0$ with
$|e_{k_0}(x)| \ge \frac 1{k_0}$  and $|y(x)|\ge \rho^{1/4}$
for $|x|\le \rho$ and $|x-\rho|\le \varepsilon$.
These assumptions imply
\[
|e_{k_0}(x)| \ge \frac 1{k_0} \ge e^{-2C_2 \sqrt {k_0}} \ge
|y(x)|^{k_0} e^{-2C_2 \sqrt {k_0}}
\]
and by Corollary~\ref{CorEk} (since $|x|\le\rho<\sqrt\rho$ this is applicable)
\[
\frac{|E_{k_0}|}{|e_{k_0}^2|} \le C_0 \rho^{k_0} |y|^{-2k_0} e^{4C_2 \sqrt k_0} \le
C_0 \rho^{k_0/2} e^{4C_2 \sqrt k_0} \le 1.
\]

The goal is to show by induction that for $k\ge k_0$
and for $|x|\le \rho$ and $|x-\rho|\le \varepsilon$
\begin{equation}\label{eqekgoal}
|e_k| \ge |y|^k e^{- 3C_2 \sqrt k} \quad\mbox{and}\quad
\left| \frac{E_k}{e_k^2}\right| \le 1.
\end{equation}
Assume that (\ref{eqekgoal}) is satisfied for $k= k_0$.
Now suppose that (\ref{eqekgoal}) holds for some $k\ge k_0$.
Then (\ref{eqek1ek}) implies
\begin{align*}
|e_{k+1}|&\ge |y|\, |e_k| e^{-C_2 k^{-1/2}}  \\
&\ge |y|^{k+1} e^{- 3C_2 \sqrt k} e^{-C_2 k^{-1/2}}  \\
&\ge |y|^{k+1} e^{- 3C_2 \sqrt {k+1}}.
\end{align*} 
Furthermore
\[
\frac{|E_{k+1}|}{|e_{k+1}^2|} \le C_0 \rho^{k+1} |y|^{-2k-2} e^{4C_2 \sqrt {k+1}} \le
C_0 \rho^{(k+1)/2} e^{4C_2 \sqrt {k+1}} \le 1.
\]
Hence, we have proved (\ref{eqekgoal}) for all $k\ge k_0$.

In the last step of the induction proof we also obtained the 
upper bound
\[
\frac{|E_k|}{|e_k^2|}\le  C_0 \rho^{k/2} e^{4C_2 \sqrt k} 
\]
which is sufficient to obtain the asymptotic representation
(\ref{eqekrep}).
\epf

\bl\label{Leek2}
Suppose that $|x|\le \rho-\varepsilon$ for some $\varepsilon> 0$.
Then we have uniformly
\begin{equation} \label{ekdarstellung}
e_k(x) = C_k(x) y(x)^k = (C(x) + o(1))y(x)^k
\end{equation} 
for some analytic function $C(x)$. Consequently we have
uniformly for $|x|\le \sqrt{\rho-\varepsilon}$
\begin{equation} \label{Ekdarstellung} 
E_k(x) = \tilde C_k(x) y(x^2)^k = (\tilde C(x) + o(1))y(x^2)^k
\end{equation} 
with an analytic function $\tilde C(x)$.
\el

\bpf
If $|x|\le \rho - \varepsilon$ then by Lemma~\ref{upperbound_height} we have
$|e_k(x)|\le e_k(\rho-\varepsilon) = \Ord{\(1-\frac\eps\rho\)^k}$.
Thus, we can replace the upper bound $|e_k(x)| \le C/\sqrt k$
in the proof of Lemma~\ref{Leek1} by an exponential bound
which leads to a lower bound for $e_k(x)$ of the form
\[
|e_k(x)|\ge c_0 |y(x)|^k.
\]
Hence, by using (\ref{eqekexprep2}) the result follows
with straightforward calculations.

In order to show the second assertion we start from \eqref{Ekdef} and insert 
\eqref{ekdarstellung}. Then we obtain 
\begin{align} 
E_k(x)&= \frac{C(x^2) + o(1)}2 y(x^2)^k+\sum_{i\ge 3}\frac{(C(x^i) + o(1))y(x^i)^k}i \nonumber \\
&=\(\frac{C(x^2)}2+ \sum_{i\ge 3}\frac{(C(x^i) + o(1))y(x^i)^k}{iy(x^2)} + o(1)\) y(x^2).
\label{Ekdetail}
\end{align}
Since $C(x)$ is analytic, it is bounded in the compact interval $[0,\rho-\eps]$. Furthermore,
observe that using the bound from Lemma~\ref{lem_2} it is easy to see that 
the series in \eqref{Ekdetail} is uniformly convergent. Hence the representation
\eqref{Ekdarstellung} follows.  
\epf

The disadvantage of the previous two lemmas is that they only work
for $|x|\le \rho$.
In order to obtain some progress for $|x|>\rho$ fix a constant $C>0$ such
that 
\[
\left| e^{-E_k(x)} -1 \right| \le 
\frac C{\sqrt k} \left( \frac{|x|^2}{\rho} \right)^k
\]
for all $k\ge 1$ and for all $|x|\le \sqrt \rho$.

\bl\label{LeGWh3.2}
Let $x\in \Delta$ and suppose that there exist real numbers
$D_1$ and $D_2$ with $0< D_1,D_2 < 1$ and some integer $K\ge 1$ with
\begin{equation}\label{eqLeGWh3.2}
|e_{K}(x)| < D_1, \quad |y(x)| \frac{e^{D_1} - 1}{D_1} < D_2, \quad
D_1 D_2 + e^{D_1} \frac C{\sqrt {K}} 
\left( \frac{|x|^2}{\rho} \right)^{K} < D_1.
\end{equation}
Then we have $|e_k(x)|< D_1$ for all $k\ge K$ and
\[
e_k(x) =  O(y(x)^k)
\]
as $k\to\infty$, where the implicit constant might depend on $x$.
\el

\brem
Note that the second inequality in \eqref{eqLeGWh3.2} implies $|y(x)|<1$ and thus an exponential
decay of $e_k(x)$. Hence, in particular, the assumptions of this lemma imply that
\eqref{einsdurchk} holds.
\erem

\bpf
By definition we have $e_{K+1} = y(1 - e^{-e_{K} - E_{K}})$. 
Hence, if we write $e^{-E_{k}} = 1 + R_{k}$ we  obtain
\[
|e_{K+1}| \le |e_{K}|\, |y|\frac{e^{|e_{K}|} - 1}{|e_{K}|}
 + e^{|e_{K}|} R_{K}.
\]
If (\ref{eqLeGWh3.2}) is satisfied then it follows that
\[
|e_{K+1}| \le D_1 D_2 + e^{D_1} \frac C{\sqrt {K}} 
\left( \frac{|x|^2}{\rho} \right)^{K} < D_1.
\]
Now we can proceed by induction and obtain $|e_k|< D_1$ for
all $k\ge K$. Note that $D_2 < 1$ and since $x\in \Delta$ we have $|x|\le \rho+\eps$,
Corollary~\ref{CorEk} implies 
$$
|E_k(x)|\le \frac{C_0}{\sqrt k} \left|\frac{(\rho+\eps)^2}{\rho}\right|^k<(\rho+3\eps)^k.
$$
Moreover, $R_k=e^{-E_k}-1=\Ord{|E_k|}$ and hence $\sum_k R_k = O(1)$. Thus we have 
\[
e_k = O(D_2^k).
\]

If we set $a_k = y (1-e^{-e_k})/e_k$ and 
$b_k = - ye^{-e_k}R_k$ we obtain the recurrence
\[
e_{k+1} = e_k a_k + b_k
\]
with an explicit solution of the form
\[
e_k = e_{K}\prod_{K\le i < k} a_i + \sum_{K\le j < k} b_j 
\prod_{j < i < k} a_i.
\]
Since $e_k = O(D_2^k)$, we have 
\[
\prod_{j < i < k} a_i = \Ord{y(x)^{k-j}}
\] 
and hence by Corollary~\ref{CorEk} 
\[
b_j  \prod_{j < i < k} a_i = O\left( y(x)^{k-j} |x^2/\rho|^j \right) =
O\left( y(x)^k L^j \right)
\]
for some $L$ with $0< L < 1$. Hence,
\begin{equation}\label{eqekasyrep}
e_k = \Ord{e_{K} y(x)^{k-K}} + O\left( y(x)^k L^{K} \right) 
= O\left(y(x)^k\right).
\end{equation}
\epf

Recall that $e_k(x)\to 0$ for $|x|\le \rho$.
Using Lemma~\ref{LeGWh3.2} we deduce that
$e_k\to 0$ in a certain region that extends the circle $|x|\le \rho$.
\bl
For every $\varepsilon>0$ there exists $\delta>0$ such that
$e_k(x)\to 0$ at an exponential rate, if $|x|\le \rho+\delta$ and $|x-\rho|\ge \varepsilon$.
\el
\bpf
If we show that Lemma~\ref{LeGWh3.2} is applicable, then we are done. 
First observe that since the assumption on $x$ implies $|y(x)|<1$, for 
sufficiently small $D_1$ there exist $D_2<1$ and $K$ (sufficiently large) such that the 
second and third inequality if \eqref{eqLeGWh3.2} hold. Next note that $|e_K(x)|\le 
e_K(\rho)$ for $|x|\le \rho$ and thus $|e_K(x)|\le 2e_K(\rho)$ by continuity. By \eqref{ek2k} we
have $e_K(\rho)\sim\frac 2K$ we can make $|e_K(x)|$ arbitrarily small and the proof is complete. 
\epf

Now we turn to the most important range, namely for $x\in \Delta$ 
and $|x-\rho|\le \varepsilon$. 

\bl\label{LeGWh6.2}
There exists $\varepsilon>0$ and a constant $c_1>0$ such that
for all $x\in \Delta_\eps$ 
the conditions (\ref{eqLeGWh3.2}) are satisfied 
for  
\[
k =K(x) = \left\lfloor \frac {c_1}{|\arg (y(x))|}\right\rfloor
\]
and properly chosen real numbers $D_1,D_2$. Consequently,
$e_k(x)\to 0$ at exponential rate. 
\el

\bpf
Suppose that $\arg(x)$ and $\arg(y(x))$ are positive
and that $K(x) = $ \linebreak $\left\lfloor  {c_1}/{|\arg (y(x))|}\right\rfloor$,
where $c_1 = \arccos(1/4)- \varepsilon_1$ and 
$\varepsilon_1>0$ is arbitrarily small. Note that $K(x)$ can be made as large as we desire, since
$y(x)\sim 1$ and therefore $\arg(y(x))$ is small for small $\eps$.

Now fix an integer $k_0$ and $\eps$ small enough to guarantee $k_0<K(x)$.
Moreover, fix two small positive real numbers $\eps_2$ and $\eps_3$. First we
will prove by induction that for $k_0\le k\le K(x)$ we have (for sufficiently small $\eps$)
\begin{align} 
|e_k(x)|&\le \eps_2, \label{f1} \\
\arg(e_k(x))&\le k\arg(y(x))+\eps_3, \label{f2} \\
|e_k(x)|&\le c\frac{|y(x)|^k}k \mbox{ for some $c>0$.} \label{f3}
\end{align}

The first step of the induction proof is to show \eqref{f1}. Observe that due to the choice of
$c_1$ and formula \eqref{f2} of the induction hypothesis we have 
\begin{equation} 
0<\arg(e_k(x)) < \arccos(1/4). \label{arg}
\end{equation} 
This implies $|e_{k+1}| \le |e_k| + |E_k|$ and consequently, by
the second statement of 
Lemma~\ref{Leek2} and the property $e_{k_0}(\rho) \sim c/k_0$ (compare with \eqref{ek2k})
it follows that
\begin{equation} \label{auxi}
|e_k(x)| \le |e_{k_0}| + \sum_{k_0 \le \ell < k} |E_\ell| < \varepsilon_2
\end{equation}
provided that $k_0$ is chosen sufficiently large.

Next we show \eqref{f2}. We start with \eqref{eqrecek} and obtain 
\begin{align} 
e_{k+1}&=y(x) e_k\(\frac{1-e^{-e_k}}{e_k}+\Ord{\frac{E_k}{e_k}}\) \nonumber \\
&=y(x) e_k\(1-\frac{e_k}2+\Ord{e_k^2} +\Ord{\frac{E_k}{e_k}}\). \label{zwischenschritt}
\end{align} 
Note that by \eqref{f1} the first of the two error terms is much smaller than $e_k$. In order to
estimate the second error term, note that by the second statement of
Lemma~\ref{Leek2} we know $E_k(x) = (\tilde C(x)+o(1)) y(x^2)^k$. Combining this with \eqref{f1}
we obtain 
$$
\left|\frac{E_k}{e_k}\right|=\Ord{\left|\frac{ky(x^2)^k}{y(x)^k}\right|}=\Ord{L^k} 
$$
with some $0<L<1$. Since $|y(x^2)/y(x)|<\rho$ (a consequence of the convexity of $y(x)$ on the
positive real line) whereas $|y(x)|\sim \rho$, we can have $L<|y(x)|$ provided that $\eps$ is
small enough. But this together with \eqref{f3} implies that also the second error term in
\eqref{zwischenschritt} is small in comparison to $e_k$. Hence \eqref{arg} implies that the
argument of the last factor in \eqref{zwischenschritt} is negative. Thus we conclude 
\begin{align*} 
\arg(e_{k+1}(x)) &\le \arg(y(x))+ \arg(e_k) \\
&\le (k-k_0+1)\arg(y(x)) + \arg(e_{k_0}(x)) \\
&\le k\arg(y(x)) + \varepsilon_3 
\end{align*} 
where the last inequality follows the fact that by 
$e_{k_0}(\rho) \sim c/k_0$ (compare with \eqref{ek2k}) and continuity we can always achieve
$|\arg(e_{k_0}(x))|<\eps_3$. 

The third step is to prove the lower bound \eqref{f3} for $e_k(x)$ for $k_0\le k\le K(x)$.
of the form $|e_{k}(x)|\ge c |y(x)|^{k}/k$ for some $c>0$.
By Lemma~\ref{Leek2} $E_k(x) = (\tilde C(x)+o(1))
y(x^2)^k$ behaves nicely, if $|x-\rho|\le \varepsilon$. 
Suppose that $|x|\ge \rho$ and $x\in \Delta_\eps$.
Since $\arg(y(x^2))$ is of order $\arg(y(x))^2$ we deduce that
$\arg(E_k(x)) = O(\arg(y(x)))$ for $k_0\le k \le K(x)$. 
In particular, it follows that (for $k_0\le k \le K(x)$)
\[
|e_{k+1}(x)| = \left|y(x)(1-e^{-e_k(x)-E_k(x)})\right| \ge  
|y(x)|(1-e^{-|e_k(x)|}).
\]
Treating the nonlinear recurrence $a_{k+1}=1-e^{-a_k}$ 
with the methods of de Bruijn \cite[p. 156]{dB}, it is possible to show inductively that 
$a_k\sim \frac ck$ and thus $|e_{k}(x)|\ge c |y(x)|^{k}/k$ for some $c>0$. 

The last task is to find $D_1$ and $D_2$ such that the conditions (\ref{eqLeGWh3.2}) 
are satisfied for $k = K(x)$. In order to do this, we first show that in the formula 
\begin{equation}\label{eqekexprep2-2}
e_k = \frac{y^{k-k_0}}{\frac 1{e_{k_0}} + \frac 12 \frac{1-y^{k-k_0}}{1-y}
+ \sum_{k_0\le \ell < k} e_\ell h(e_\ell) y^\ell + 
\sum_{k_0\le \ell < k} \frac{S_\ell}{e_\ell^2} y^\ell  }, 
\end{equation}
where $k_0$ is fixed, the second term in the denominator dominates.
Since $k_0$ is fixed, the first term is bounded. For estimating the third term note that by
\eqref{auxi} the terms in the sum can be made arbitrarily small if $k_0$ is chosen sufficiently
large. Finally, due to the already obtained bounds $|e_{k}(x)|\ge c |y(x)|^{k}/k$,
$E_k(x) = O(y(x^2)^k)$, and the property $S_\ell \sim - E_\ell$
the last term satisfies
\[
\sum_{k_0\le \ell < k} \frac{S_\ell}{e_\ell^2} y^\ell = O(1)
\]
and therefore it does not contribute to the main term, either. Summing up we have
\[
|e_k| \le \frac{|y|^{k-k_0}}{\left| \frac 12 \frac{1-y^{k-k_0}}{1-y}  \right|} (1 + \varepsilon_4)
\]
for an arbitrarily small $\varepsilon_4>0$.

We set $\arg(\rho-x) = \theta$, where we assume that
$\theta\in \left[-\frac \pi 2 - \varepsilon_5,
\frac \pi 2 + \varepsilon_5 \right]$ (for some $\varepsilon_5> 0$ that
has to be sufficiently small), and $r = b|\rho-x|^{1/2}$, where
$b$ is the constant appearing in \eqref{y-asym}.
Then we have 
\begin{align*}
|y| &= 1 - r \cos \frac \theta 2 + O(r^2),\\
\log |y| &= -  r \cos \frac \theta 2 + O(r^2),\\
\arg(y) &=  -r \sin \frac \theta 2 + O(r^2).
\end{align*}

Hence with $k = K(x)= \lfloor c_1/|\arg(y)|\rfloor$ we have
\[
|y^{k-k_0}| \sim e^{-c_1 \cot(\theta/2) + O(r^2)} 
\le e^{-c_1 \cot(\frac \pi 4 + \frac{\varepsilon_5}2) + O(r^2)} \le
e^{-c_1}(1-\varepsilon_6)
\]
for some arbitrarily small $\varepsilon_6>0$ (depending on
$\varepsilon_5$). Consequently
\[
|e_k| < D_1 := 2 \frac{e^{-\arccos(1/4)}}{1- e^{-\arccos(1/4)}} r 
(1 +  \varepsilon_7) = c' r,
\]
where $\varepsilon_7>0$ can be chosen arbitrarily small.
Moreover 
\[
|y| = 1 - r \cos \frac \theta 2 + O(r^2) \le 1 - \frac r{\sqrt 2}
(1- \varepsilon_8)
\]
for some (small) $\varepsilon_8>0$ and consequently
\[
|y| \frac{e^{D_1} - 1}{D_1} = 1 - \left( \frac 1{\sqrt 2} - 
\frac{e^{-\arccos(1/4)}}{1- e^{-\arccos(1/4)}} \right) r (1-\varepsilon_9) 
+ O(r^2).
\]
Thus, we are led to set
\[
D_2 := 1 - \frac 12\left( \frac 1{\sqrt 2} - 
\frac{e^{-\arccos(1/4)}}{1- e^{-\arccos(1/4)}} \right) r = 1 - c'' r
\]
and the first two conditions of (\ref{eqLeGWh3.2}) are satisfied
if $r$ is sufficiently small.
Since $D_1 - D_1D_2 = c'c'' r^2$ we just have to check whether
\[
e^{D_1} \frac C{\sqrt k} \left| \frac{x^2}\rho \right|^k < c'c'' r^2.
\]
However, since $k = K(x) \ge c_1 \sqrt 2\, r^{-1}$ the left hand side
of this inequality 
is definitely smaller than $c'c'' r^2$ if $r$ is sufficiently small. 
Hence all conditions of (\ref{eqLeGWh3.2}) are satisfied for
$k = K(x)$.
\epf

\bl\label{Leek1-2}
There exists $\varepsilon>0$ such that (\ref{eqekrep})
holds for all $x$ with $x\in \Delta_\eps$ with 
$|x|\ge \rho$.
\el

\bpf
We recall that the properties $e_k = O(1/k)$ and (\ref{eqProOLk}) imply
(\ref{eqekrep}). By Lemma~\ref{LeGWh6.2} we already know that
the first condition holds. Furthermore, we have upper bounds for $E_k$ (see Lemma~\ref{Leek2}).
Hence, it remains to provide proper lower bounds for $e_k$.

Since we already know that $e_k\to 0$ and $|E_k| = O(L^k)$
(for some $L<1$), the
recurrence (\ref{eqrecek}) implies
\[
|e_{k+1}| \ge (1-\delta)\left( |e_k| - |E_k|\right)
\]
for some $\delta>0$ provided that $x\in \Delta$ and $|x-\rho|< \varepsilon$.
Without loss of generality we can assume that $L < (1-\delta)^2$.
Hence
\[
|e_k| \ge (1-\delta)^k - \sum_{\ell < k} |E_\ell| (1-\delta)^{k-\ell}
\ge c_0 (1-\delta)^k
\]
for some constant $c_0>0$. Consequently
\[
\left| \frac{E_k}{e_k^2} \right| = O\left( \left( \frac L{(1-\delta)^2} \right)^k
\right).
\]
As noted above, this upper bound is sufficient to deduce (\ref{eqekrep}).
\epf

\begin{ack} 
We express our gratitude to Nicolas Broutin for his careful reading of an earlier version of our
manuscript as well as to him and Philippe Flajolet for inspiring discussions on the
topic and for providing us with a preprint of their paper \cite{broutinflajolet}. Furthermore, we
thank two anonymous referees whose detailed and very informative reports led to the correction 
of numerous imprecisions and errors in the manuscript as well as to a considerable 
improvement of the presentation. 
\end{ack}


\begin{thebibliography}{10}

\bibitem{Al91}
David Aldous.
\newblock The continuum random tree. {II}. {A}n overview.
\newblock In {\em Stochastic analysis (Durham, 1990)}, volume 167 of {\em
  London Math. Soc. Lecture Note Ser.}, pages 23--70. Cambridge Univ. Press,
  Cambridge, 1991.

\bibitem{Bi68}
Patrick Billingsley.
\newblock {\em Convergence of probability measures}.
\newblock John Wiley \& Sons Inc., New York, 1968.

\bibitem{broutinflajolet}
N.~Broutin and P.~Flajolet.
\newblock The height of random binary unlabelled trees.
\newblock {\em Discrete Math. Theor. Comput. Sci. Proc.}, AI:121--134, 2008.

\bibitem{CMY}
P.~Chassaing, J.~F. Marckert, and M.~Yor.
\newblock The height and width of simple trees.
\newblock In {\em Mathematics and computer science (Versailles, 2000)}, Trends
  Math., pages 17--30. Birkh\"auser, Basel, 2000.

\bibitem{CDJ01}
Brigitte Chauvin, Michael Drmota, and Jean Jabbour-Hattab.
\newblock The profile of binary search trees.
\newblock {\em Ann. Appl. Probab.}, 11(4):1042--1062, 2001.

\bibitem{CoHo81}
J.~W. Cohen and G.~Hooghiemstra.
\newblock Brownian excursion, the {$M/M/1$} queue and their occupation times.
\newblock {\em Math. Oper. Res.}, 6(4):608--629, 1981.

\bibitem{dB}
N.~G. de~Bruijn.
\newblock {\em Asymptotic methods in analysis}.
\newblock Dover Publications Inc., New York, third edition, 1981.

\bibitem{BKR}
N.~G. de~Bruijn, D.~E. Knuth, and S.~O. Rice.
\newblock The average height of planted plane trees.
\newblock In {\em Graph theory and computing}, pages 15--22. Academic Press,
  New York, 1972.

\bibitem{DeHw06}
Luc Devroye and Hsien-Kuei Hwang.
\newblock Width and mode of the profile for some random trees of logarithmic
  height.
\newblock {\em Ann. Appl. Probab.}, 16(2):886--918, 2006.

\bibitem{Dr96}
M.~Drmota.
\newblock On nodes of given degree in random trees.
\newblock In {\em Probabilistic methods in discrete mathematics (Petrozavodsk,
  1996)}, pages 31--44. VSP, Utrecht, 1997.

\bibitem{Dr04}
Michael Drmota.
\newblock On {R}obson's convergence and boundedness conjectures concerning the
  height of binary search trees.
\newblock {\em Theoret. Comput. Sci.}, 329(1-3):47--70, 2004.

\bibitem{Drm09}
Michael Drmota.
\newblock {\em Random trees}.
\newblock SpringerWienNewYork, Vienna, 2009.
\newblock An interplay between combinatorics and probability.

\bibitem{DG97}
Michael Drmota and Bernhard Gittenberger.
\newblock On the profile of random trees.
\newblock {\em Random Structures Algorithms}, 10(4):421--451, 1997.

\bibitem{DG3}
Michael Drmota and Bernhard Gittenberger.
\newblock The distribution of nodes of given degree in random trees.
\newblock {\em J. Graph Theory}, 31(3):227--253, 1999.

\bibitem{DG04}
Michael Drmota and Bernhard Gittenberger.
\newblock The width of {G}alton-{W}atson trees conditioned by the size.
\newblock {\em Discrete Math. Theor. Comput. Sci.}, 6(2):387--400 (electronic),
  2004.

\bibitem{DH05a}
Michael Drmota and Hsien-Kuei Hwang.
\newblock Bimodality and phase transitions in the profile variance of random
  binary search trees.
\newblock {\em SIAM J. Discrete Math.}, 19(1):19--45 (electronic), 2005.

\bibitem{DH05b}
Michael Drmota and Hsien-Kuei Hwang.
\newblock Profiles of random trees: correlation and width of random recursive
  trees and binary search trees.
\newblock {\em Adv. in Appl. Probab.}, 37(2):321--341, 2005.

\bibitem{FO90}
Ph. Flajolet and A.~M. Odlyzko.
\newblock Singularity analysis of generating functions.
\newblock {\em SIAM Journal on Discrete Mathematics}, 3:216--240, 1990.

\bibitem{FGOR93}
Philippe Flajolet, Zhicheng Gao, Andrew Odlyzko, and Bruce Richmond.
\newblock The distribution of heights of binary trees and other simple trees.
\newblock {\em Combin. Probab. Comput.}, 2(2):145--156, 1993.

\bibitem{FO82}
Philippe Flajolet and Andrew Odlyzko.
\newblock The average height of binary trees and other simple trees.
\newblock {\em J. Comput. System Sci.}, 25(2):171--213, 1982.

\bibitem{FlSe10}
Philippe Flajolet and Robert Sedgewick.
\newblock {\em Analytic combinatorics}.
\newblock Cambridge University Press, Cambridge, 2009.

\bibitem{FHN06}
Michael Fuchs, Hsien-Kuei Hwang, and Ralph Neininger.
\newblock Profiles of random trees: limit theorems for random recursive trees
  and binary search trees.
\newblock {\em Algorithmica}, 46(3-4):367--407, 2006.

\bibitem{G98}
Bernhard Gittenberger.
\newblock Convergence of branching processes to the local time of a {B}essel
  process.
\newblock In {\em Proceedings of the Eighth International Conference ``Random
  Structures and Algorithms'' (Poznan, 1997)}, volume~13, pages 423--438, 1998.

\bibitem{G02}
Bernhard Gittenberger.
\newblock On the profile of random forests.
\newblock In {\em Mathematics and computer science, II (Versailles, 2002)},
  Trends Math., pages 279--293. Birkh\"auser, Basel, 2002.

\bibitem{Gi06}
Bernhard Gittenberger.
\newblock Nodes of large degree in random trees and forests.
\newblock {\em Random Structures Algorithms}, 28(3):374--385, 2006.

\bibitem{GL99}
Bernhard Gittenberger and Guy Louchard.
\newblock The {B}rownian excursion multi-dimensional local time density.
\newblock {\em J. Appl. Probab.}, 36(2):350--373, 1999.

\bibitem{HRS75}
Frank Harary, Robert~W. Robinson, and Allen~J. Schwenk.
\newblock Twenty-step algorithm for determining the asymptotic number of trees
  of various species.
\newblock {\em J. Austral. Math. Soc. Ser. A}, 20(4):483--503, 1975.

\bibitem{Hoog82}
G.~Hooghiemstra.
\newblock On the explicit form of the density of {B}rownian excursion local
  time.
\newblock {\em Proc. Amer. Math. Soc.}, 84(1):127--130, 1982.

\bibitem{Hw05}
Hsien-Kuei Hwang.
\newblock Profiles of random trees: plane-oriented recursive trees (extended
  abstract).
\newblock In {\em 2005 International Conference on Analysis of Algorithms},
  Discrete Math. Theor. Comput. Sci. Proc., AD, pages 193--200 (electronic).
  Assoc. Discrete Math. Theor. Comput. Sci., Nancy, 2005.

\bibitem{Hw07}
Hsien-Kuei Hwang.
\newblock Profiles of random trees: plane-oriented recursive trees.
\newblock {\em Random Structures Algorithms}, 30(3):380--413, 2007.

\bibitem{KS}
Ioannis Karatzas and Steven~E. Shreve.
\newblock {\em Brownian motion and stochastic calculus}, volume 113 of {\em
  Graduate Texts in Mathematics}.
\newblock Springer-Verlag, New York, 1988.

\bibitem{Kol77}
V.~F. Kolchin.
\newblock Branching processes, random trees and a generalized particle
  distribution scheme.
\newblock {\em Mat. Zametki}, 21(5):691--705, 1977.

\bibitem{Kol86}
Valentin~F. Kolchin.
\newblock {\em Random mappings}.
\newblock Translation Series in Mathematics and Engineering. Optimization
  Software Inc. Publications Division, New York, 1986.

\bibitem{LST}
Guy Louchard, Wojciech Szpankowski, and Jing Tang.
\newblock Average profile of the generalized digital search tree and the
  generalized {L}empel-{Z}iv algorithm.
\newblock {\em SIAM J. Comput.}, 28(3):904--934 (electronic), 1999.

\bibitem{MaMi09}
J.-F. Marckert and G.~Miermont.
\newblock The crt is the scaling limit of unordered binary trees.

\bibitem{MM78}
A.~Meir and J.~W. Moon.
\newblock On the altitude of nodes in random trees.
\newblock {\em Canadian Journal of Mathematics}, 30:997--1015, 1978.

\bibitem{MM91}
A.~Meir and J.~W. Moon.
\newblock On nodes of large out-degree in random trees.
\newblock In {\em Proceedings of the Twenty-second Southeastern Conference on
  Combinatorics, Graph Theory, and Computing (Baton Rouge, LA, 1991)},
  volume~82, pages 3--13, 1991.

\bibitem{MM92}
A.~Meir and J.~W. Moon.
\newblock On nodes of given out-degree in random trees.
\newblock In {\em Fourth Czechoslovakian Symposium on Combinatorics, Graphs and
  Complexity (Prachatice, 1990)}, volume~51 of {\em Ann. Discrete Math.}, pages
  213--222. North-Holland, Amsterdam, 1992.

\bibitem{Ni05}
Pierre Nicod{\`e}me.
\newblock Average profiles, from tries to suffix-trees.
\newblock In {\em 2005 International Conference on Analysis of Algorithms},
  Discrete Math. Theor. Comput. Sci. Proc., AD, pages 257--266 (electronic).
  Assoc. Discrete Math. Theor. Comput. Sci., Nancy, 2005.

\bibitem{OW87}
Andrew~M. Odlyzko and Herbert~S. Wilf.
\newblock Bandwidths and profiles of trees.
\newblock {\em J. Combin. Theory Ser. B}, 42(3):348--370, 1987.

\bibitem{Otter}
Richard Otter.
\newblock The number of trees.
\newblock {\em Ann. Math.}, 49(2):583--599, 1948.

\bibitem{Pi99}
Jim Pitman.
\newblock The {SDE} solved by local times of a {B}rownian excursion or bridge
  derived from the height profile of a random tree or forest.
\newblock {\em Ann. Probab.}, 27(1):261--283, 1999.

\bibitem{Pol37}
George P\'olya.
\newblock Kombinatorische {A}nzahlbestimmungen f\"ur {G}ruppen, {G}raphen und
  chemische {V}erbindungen.
\newblock {\em Acta Math.}, 68:145--254, 1937.

\bibitem{RS75}
Robert~W. Robinson and Allen~J. Schwenk.
\newblock The distribution of degrees in a large random tree.
\newblock {\em Discrete Math.}, 12(4):359--372, 1975.

\bibitem{Stepanov69}
V.~E. Stepanov.
\newblock The distribution of the number of vertices in the layers of a random
  tree.
\newblock {\em Teor. Verojatnost. i Primenen.}, 14:64--77, 1969.

\bibitem{Ta91}
Lajos Tak{\'a}cs.
\newblock Conditional limit theorems for branching processes.
\newblock {\em J. Appl. Math. Stochastic Anal.}, 4(4):263--292, 1991.

\bibitem{Ta91b}
Lajos Tak{\'a}cs.
\newblock On the distribution of the number of vertices in layers of random
  trees.
\newblock {\em J. Appl. Math. Stochastic Anal.}, 4(3):175--186, 1991.

\bibitem{HHM02}
Remco van~der Hofstad, Gerard Hooghiemstra, and Piet Van~Mieghem.
\newblock On the covariance of the level sizes in random recursive trees.
\newblock {\em Random Structures Algorithms}, 20(4):519--539, 2002.

\end{thebibliography}

\bibliographystyle{plain}

\def\cprime{$'$}

\end{document}